\documentclass[12pt]{article}
\usepackage{amsfonts}

\usepackage{latexsym}
\usepackage{amssymb}
\usepackage{epsfig}
\usepackage{amsmath}
\usepackage{amsthm}
\usepackage[mathscr]{eucal}
\usepackage{subfigure}
\usepackage{graphicx}

\newtheorem{Lemma}{Lemma}
\newtheorem{Proposition}{Proposition}
\newtheorem{Theorem}[Lemma]{Theorem}

\newtheorem{Corollary}[Lemma]{Corollary}

\renewcommand{\qed}{\hfill{\ \ \rule{2mm}{2mm}} \vspace{0.2in}}

\newcommand{\ind}{1\hspace{-2.3mm}{1}}

\begin{document}

\title{Cliques and Chromatic Number in Inhomogenous Random Graphs}
\author{ \textbf{Ghurumuruhan Ganesan}
\thanks{E-Mail: \texttt{gganesan82@gmail.com} } \\
\ \\
New York University, Abu Dhabi.}
\date{}
\maketitle

\begin{abstract}
In this paper, we study cliques and chromatic number of inhomogenous random graphs where the individual edge probabilities could be arbitrarily low. We use a recursive method to obtain estimates on the maximum clique size under a mild positive average edge density assumption. As a Corollary, we also obtain uniform bounds on the maximum clique size and chromatic number for homogenous random graphs for all ranges of the edge probability~\(p_n\) satisfying~\(\frac{1}{n^{\alpha_1}} \leq p_n \leq 1-\frac{1}{n^{\alpha_2}}\) for some positive constants~\(\alpha_1\) and~\(\alpha_2.\)


\vspace{0.1in} \noindent \textbf{Key words:} Random graphs, inhomogenous edge probabilities, cliques, chromatic number.

\end{abstract}

\bigskip

\setcounter{equation}{0}
\renewcommand\theequation{\thesection.\arabic{equation}}
\section{Introduction}\label{intro}


Let~\(K_n\) be the labelled complete graph on~\(n\) vertices with vertex set\\\(\{1,2,\ldots,n\}\) and edge set~\(\{e_1,e_2,\ldots,e_m\},\) where~\(m = {n \choose 2}.\) Let~\(G_n = G(n,p_n)\) be the random graph obtain when every edge is independently open with probability~\(p_n \in (0,1)\) and closed otherwise. Let~\(X(i,j)\) be a Bernoulli random variable defined on the probability space\\\((\{0,1\}, \mathbb{B}(\{0,1\}), \mathbb{P}_{i,j})\) with~\[\mathbb{P}_{i,j}(X(i,j) = 1) = p(i,j) = 1 -\mathbb{P}_{i,j}(X(i,j) = 0).\] Here~\(\mathbb{B}(\{0,1\})\) is the set of all subsets of~\(\{0,1\}.\) We say that edge~\(e(i,j)\) is \emph{open} if~\(X(i,j) = 1\) and closed otherwise. The random variables~\(\{X(i,j)\}\) are independent and the resulting random graph~\(G\) is an inhomogenous Erd\H{o}s-R\'enyi (ER) random graph, defined on the probability space~\((\Omega,{\cal F}, \mathbb{P}).\) Here~\(\Omega = \{0,1\}^{ {n \choose 2}},\) the sigma algebra~\({\cal F}\) is the set of subsets of~\(\Omega\) and~\(\mathbb{P} = \prod_{i,j} \mathbb{P}_{i,j}.\)

\subsection*{Clique Number}
Let~\(G = (V,E)\) be a graph with vertex set~\(V = V(G)\) and edge set~\(E = E(G).\) Suppose that~\(\#V = n\) so that~\(G\) is a graph on~\(n\) vertices. We say that~\(G\) is a \emph{complete graph} if~\(\#E = {n \choose 2}.\) For a subset~\(U \subset V\) containing~\(\#U = r \leq n\) vertices, we define~\(G|_U = (U,E_U)\) to be the \emph{induced subgraph} of~\(G\) with vertex set~\(U,\) defined as follows. For any two vertices~\(a,b \in U,\) the edge~\(e_{ab}\) with endvertices~\(a\) and~\(b\) belongs to~\(E_U\) if and only if~\(e_{ab} \in E(G).\) We say that~\(G|_U\) is a \emph{clique} if~\(G|_U\) is a complete graph. We denote~\(\omega(G)\) to be the size of the largest clique in~\(G.\) Throughout, the size of a graph will always refer to the number of vertices in the graph.

Let~\(\textbf{p} = \{p(i,j)\}_{i,j}\) be a vector of probabilities as defined in the previous subsection and let~\(G(n,\textbf{p})\) be the resulting random graph. The following obtains an upper bound for the clique number~\(\omega(G(n,\textbf{p})).\) 
\begin{Proposition}\label{clq_up_bound}
For any sequence~\(U_n  > 0\) define
\begin{equation}\label{log_ave}
\log\left(\frac{1}{t_n}\right) := \inf_{S : \#S =U_n} \left({\#S \choose 2}\right)^{-1}\sum_{i,j \in S} \log\left(\frac{1}{p(i,j)}\right).
\end{equation}
We have
\begin{equation}\label{clq_upper}
\mathbb{P}_{\textbf{p}}\left(\omega(G(n,\textbf{p})) \leq U_n\right) \geq 1-\exp\left(-f_nU_n\right)
\end{equation}
where
\begin{equation}\label{fn_def}
f_n := \frac{(U_n-1)}{2}\log\left(\frac{1}{t_n}\right) - \log{n}
\end{equation}
for all~\(n \geq 2.\)
\end{Proposition}




To obtain a lower bound on the clique number, we have a few definitions first. As before, let~\(\textbf{p}\) be the vector formed by the probabilities~\(\{p(i,j)\}_{i,j}.\) For~\(p_n > 0\) and constant~\(0 \leq a < 1,\) let~\({\cal N}(a,n,p_n)\) be the set of all vectors~\(\textbf{p}\) satisfying the following condition: There is a constant~\(N =N(a) \geq 1\) such that for all~\(n \geq N,\) we have
\begin{equation}\label{beta_12def}
\inf_{1 \leq i \leq n}\inf_{S} \frac{1}{\#S}\sum_{j \in S} p(i,j) \geq p_n.
\end{equation}
Here~\(N = N(a) \geq 1\) is a constant not depending on~\(n.\) For a fixed~\(1 \leq i \leq n,\) the infimum above is taken over all sets~\(S\) such that~\(\#S \geq n^{a}\) and~\(i \notin S.\) The condition implies that the average edge density taken over sets of cardinality at least~\(n^{a},\) is at least~\(p_n.\) All constants mentioned throughout are independent of~\(n.\)

Let~\(p_n\) be as in~(\ref{beta_12def}) and define
\begin{equation}\label{alp1_def}
\alpha_1 = \limsup_n \frac{\log\left(\frac{1}{p_n}\right)}{\log{n}}
\end{equation}
and
\begin{equation}\label{alp2_def}
\alpha_2 = \limsup_n \frac{\log\left(\frac{1}{1-p_n}\right)}{\log{n}}.
\end{equation}
We consider three cases separately depending on whether~\(\alpha_1 > 0\) or~\(\alpha_2 > 0\) or both~\(\alpha_1 = \alpha_2 = 0.\)
\begin{Theorem}\label{clq_main_thm}
Suppose~\(\textbf{p} \in {\cal N}(a,n,p_n)\) for some constant~\(0 \leq a <1.\)\\
\((i)\) Suppose~\(0 < \alpha_1  < 2\) and let~\(\eta,\gamma >0\) be such that
\begin{equation}\label{clq_condi}
\max\left(\frac{\alpha_1}{2},a\right)  + \gamma < \eta < 1.
\end{equation}
We have that~\(\alpha_2 = 0\) and there is a positive integer~\(N_1 = N_1(\eta,\gamma,\alpha_1,a) \geq 1\) so that
\begin{equation}\label{clq_loweri}
\mathbb{P}_{\textbf{p}}\left(\omega(G(n,\textbf{p})) \geq (1-\eta) \frac{\log{n}}{\log\left(\frac{1}{p_n}\right)}\right) \geq 1-3\exp\left(-n^{2\eta - 2\gamma-\alpha_1}\right)
\end{equation}
for all~\(n \geq N_1.\) \\
\((ii)\) Suppose~\(\alpha_1=  \alpha_2 = 0.\) Let~\(\eta,\gamma >0\) be such that
\begin{equation}\label{clq_condii}
a < \gamma < \eta <1.
\end{equation}
There is a positive integer~\(N_2 = N_2(\eta,\gamma) \geq 1\) so that
\begin{equation}\label{clq_lowerii}
\mathbb{P}_{\textbf{p}}\left(\omega(G(n,\textbf{p}) \geq (1-\eta) \frac{\log{n}}{\log\left(\frac{1}{p_n}\right)}\right) \geq 1-3\exp\left(-n^{2\eta - 2\gamma}\right)
\end{equation}
for all~\(n \geq N_2.\)\\
\((iii)\) Suppose~\(0 < \alpha_2 < 1\) and let~\(\eta,\gamma >0\) be such that
\begin{equation}\label{clq_condiii}
\max\left(\gamma-  \frac{\alpha_2}{2},a\right) < \eta <1-\alpha_2.
\end{equation}
We have that~\(\alpha_1 = 0\) and there is a positive integer~\(N_3 = N_3(\eta,\gamma) \geq 1\) so that
\begin{equation}\label{clq_loweriii}
\mathbb{P}_{\textbf{p}}\left(\omega(G(n,\textbf{p})) \geq (1-\alpha_2-\eta) \frac{\log{n}}{\log\left(\frac{1}{p_n}\right)}\right) \geq 1-3\exp\left(-n^{2\eta - 2\gamma+\alpha_2}\right)
\end{equation}
for all~\(n \geq N_3.\)
\end{Theorem}
The usual method for studying the lower bound for clique numbers of homogenous random graphs uses a combination of second moment method and martingale inequalities (see for e.g., Alon and Spencer~(2003), Bollobas~(2001)). For inhomogenous graphs where the edge probabilities could be arbitrarily low, the above method is not directly applicable. We use a recursive method to obtain lower bounds on the clique numbers (see Lemma~\ref{cliq_rec_lem}, Section~\ref{pf1}).


As a consequence of our main Theorem above, we also obtain results for homogenous random graphs where the edge probabilities~\(p(i,j) = p_n\) for~\(1 \leq i \neq j \leq n.\)
\begin{Proposition}\label{clq_extrem_lem} If~\(\alpha_1 > 2,\) then fix~\(\epsilon > 0\) small so that~\(\alpha_1 -\epsilon > 2.\) We then have
\begin{equation}\label{clq_extrem1}
\mathbb{P}\left(\omega(G(n,p_n)) = 1 \right) \geq 1-\frac{1}{n^{\alpha_1-\epsilon-2}}
\end{equation}
for all~\(n\) large. Similarly, if~\(\alpha_2 > 2,\) then fix~\(\epsilon  >0\) small so that~\(\alpha_2 -\epsilon > 2.\) We then have
\begin{equation}\label{clq_extrem2}
\mathbb{P}\left(\omega(G(n,p_n)) = n \right) \geq 1-\frac{1}{n^{\alpha_2-\epsilon-2}}
\end{equation}
for all~\(n\) large. If~\(1 < \alpha_2 < 2,\) then fix~\(\epsilon > 0\) small so that~\(0 < 2-\alpha_2-2\epsilon < 2-\alpha_2+2\epsilon < 1.\)
We then have
\begin{equation}\label{clq_extrem3}
\mathbb{P}\left(\omega(G(n,p_n) \geq n-n^{2-\alpha_2+2\epsilon}\right) \geq 1-\exp\left(-n^{2-\alpha_2-2\epsilon}\right)
\end{equation}
for all~\(n\) large.

Let~\(f_n \rightarrow \infty\) be any sequence and let
\begin{equation}\label{un_def}
U_n = \frac{2\log{n} + 2f_n}{\log\left(\frac{1}{p_n}\right)} + 1 \geq 1.
\end{equation}
We have
\begin{equation}\label{clq_upper}
\mathbb{P}\left(\omega(G(n,p_n)) \leq U_n\right) \geq 1-\exp\left(-f_nU_n\right)
\end{equation}
for all~\(n \geq 2.\)
\end{Proposition}

We have the following result regarding the clique number for the cases where~\(\alpha_1 < 1\) and~\(\alpha_2 < 1.\)
\begin{Corollary}\label{clq_cor}
\((i)\) Suppose~\(p_n = \frac{1}{n^{\theta_1}}\) for some~\(0 < \theta_1 < 1.\) Fix~\(\eta,\gamma>0\) such that \(\frac{\theta_1}{2} + \gamma < \eta<1\) and~\(\xi > 0.\)
There is a positive integer~\(N_1 = N_1(\eta,\gamma,\xi) \geq 1\) so that
\begin{equation}\label{clq_bd_case1}
\mathbb{P}_{\textbf{p}}\left(\frac{1-\eta}{\theta_1} \leq \omega(G(n,p_n)) \leq \frac{(2+\xi)}{\theta_1} + 1\right) \geq 1-3\exp\left(-n^{2\eta - 2\gamma-\theta_1}\right) - n^{-\frac{\xi(2+\xi)}{\theta_1}}
\end{equation}
for all~\(n \geq N_1.\) \\
\((ii)\) Suppose~\(p_n = p \in (0,1)\) for all~\(n.\) Fix~\(0 < \gamma < \eta < 1\) and~\(\xi > 0.\) There is a positive integer~\(N_2 = N_2(\eta,\gamma,\xi) \geq 1\) so that
\begin{eqnarray}
&&\mathbb{P}\left(\frac{(1-\eta)\log{n}}{\log\left(\frac{1}{p}\right)} \leq \omega(G(n,p_n)) \leq \frac{(2+\xi)\log{n}}{\log\left(\frac{1}{p}\right)} \right) \nonumber\\
&&\;\;\;\;\geq 1-3\exp\left(-n^{2\eta - 2\gamma}\right) - \exp\left(\frac{-\xi(1+\xi)}{\log\left(\frac{1}{p}\right)}(\log{n})^2\right)\label{clq_bd_case2}
\end{eqnarray}
for all~\(N \geq N_2.\)\\
\((iii)\) Suppose~\(p_n = 1-\frac{1}{n^{\theta_2}} \) for some~\(0 < \theta_2 < 1.\) Fix~\(\eta,\gamma> 0\) so that \(\gamma -\frac{\theta_2}{2} < \eta < 1-\theta_2\) and fix~\(\xi > 0.\) There is a positive integer~\(N_3= N_3(\eta,\gamma,\xi) \geq 1\) so that
\begin{eqnarray}
&&\mathbb{P}\left((1-\theta_2 - \eta)n^{\theta_2} \log{n} \leq \omega(G(n,p_n)) \leq (2+\xi)n^{\theta_2}\log{n}\right) \nonumber\\
&&\;\;\;\;\geq 1-3\exp\left(-n^{2\eta - 2\gamma+\theta_2}\right) - \exp\left(-\xi(1+\xi)n^{\theta_2}(\log{n})^2\right)\label{clq_bd_case3}
\end{eqnarray}
for all~\(n \geq N_3.\)
\end{Corollary}

\subsection*{Chromatic Number}


We have the following result regarding the chromatic number for homogenous random graphs where each edge is independently open with
probability~\(r_n.\) We discuss separate cases depending on the asymptotic behaviour of~\(r_n.\)
\begin{Theorem}\label{chr_cor} 
\((i)\) Suppose~\(r_n = \frac{1}{n^{\theta_2}}\) for some~\(0 < \theta_2 < \frac{1}{2}.\) Fix~\(\xi,\zeta > 0.\)
There is a constant~\(N_1 = N_1(\xi,\theta_2) \geq 1\) so that
\begin{eqnarray}&&\mathbb{P}\left((1-\xi)\frac{n^{1-\theta_2}}{2\log{n}} \leq \chi(G(n,r_n)) \leq \frac{2(1+\xi)}{1-2\theta_2} \frac{n^{1-\theta_2}}{\log{n}} \right) \nonumber\\
&&\;\;\;\;\geq 1-3\exp\left(-n^{1-\theta_2-\zeta}\right) - \exp\left(-\xi(1+\xi)n^{\theta_2}(\log{n})^2\right)\label{chr_bd_casei}
\end{eqnarray}
for all~\(n \geq N_1.\)\\
\((ii)\) Suppose~\(r_n = p\) for some~\(0 < p < 1\) and for all~\(n.\) Fix~\(\xi,\zeta > 0.\)
There is a constant~\(N_2 = N_2(\xi,\zeta) \geq 1\) so that
\begin{eqnarray}&&\mathbb{P}\left((1-\xi)\frac{n\log\left(\frac{1}{1-p}\right)}{2\log{n}} \leq \chi(G(n,r_n)) \leq 2(1+\xi) \frac{n\log\left(\frac{1}{1-p}\right)}{\log{n}} \right) \nonumber\\
&&\;\;\;\;\geq 1-3\exp\left(-n^{1 - \zeta}\right) - \exp\left(-\frac{\xi(1+\xi)}{\log\left(\frac{1}{1-p}\right)}(\log{n})^2\right)\label{chr_bd_caseii}
\end{eqnarray}
for all~\(n \geq N_2.\)\\
\((iii)\) Suppose~\(r_n = 1-\frac{1}{n^{\theta_1}}\) for some~\(0 < \theta_1 < 1.\) Fix~\(\xi,\zeta > 0.\)
There is a constant~\(N_3 = N_3(\xi,\zeta) \geq 1\) so that
\begin{eqnarray}&&\mathbb{P}\left((1-\xi)\frac{\theta_1 n}{2+\theta_1}\leq \chi(G(n,r_n)) \leq (1+\xi) \frac{2\theta_1n}{1-\theta_1} \right) \nonumber\\
&&\;\;\;\;\geq 1-3\exp\left(-n^{2\eta-2\gamma-\theta_1}\right) - n^{-\frac{\xi(1+\xi)}{\theta_1}}\label{chr_bd_caseiii}
\end{eqnarray}
for all~\(n \geq N_3.\)
\end{Theorem}


The paper is organized as follows. In Section~\ref{pf1}, we prove Proposition~\ref{clq_up_bound} and obtain preliminary estimates for proving the main Theorem~\ref{clq_main_thm}. In Section~\ref{pf_clq_main_thm}, we prove Theorem~\ref{clq_main_thm} regarding the lower bound for clique numbers of inhomogenous graphs. In Section~\ref{pf_clq_cor}, we prove Proposition~\ref{clq_extrem_lem} and Theorem~\ref{clq_cor} for clique numbers of homogenous graphs. Finally in Section~\ref{pf_chr_cor}, we prove Theorem~\ref{chr_cor} regarding the chromatic number for homogenous graphs.

\setcounter{equation}{0}
\renewcommand\theequation{\thesection.\arabic{equation}}
\section{Preliminary estimates}\label{pf1}
For integer~\(q \geq 1,\) let~\(G(q,\textbf{p})\) be the random graph with vertex set~\(S_q = \{1,2,\ldots,q\}.\)
For integer~\(L \geq 2,\) let~\(B_L(S_q)\) denote the event that the random graph~\(G(q,\textbf{p})\) contains an open~\(L-\)clique; i.e., there are vertices~\(\{v_i\}_{1 \leq i \leq L}\) such that the edge between~\(v_i\) and~\(v_j\) is open for any~\(1 \leq i \neq j \leq L.\)

\emph{Proof of Proposition~\ref{clq_up_bound}}:
We have
\begin{equation}\label{blq_upper_22}
\mathbb{P}_{\textbf{p}}(B_L(S_n)) \leq \sum_{S: \#S = L}\prod_{i,j \in S} p(i,j)  = \sum_{S: \#S = L}\exp\left(-\sum_{i,j \in S}\log\left(\frac{1}{p(i,j)}\right)\right).
\end{equation}
Setting~\(L = U_n\) and using the definition of~\(t_n\) in~(\ref{un_def}), we have
\begin{eqnarray}
\mathbb{P}_{\textbf{p}}(B_L(S_q)) &\leq& \sum_{S: \#S = U_n} \exp\left(-{U_n \choose 2} \log\left(\frac{1}{t_n}\right)\right) \nonumber\\
&\leq& {n \choose U_n} \exp\left(-{U_n \choose 2} \log\left(\frac{1}{t_n}\right)\right) \nonumber\\
&\leq& n^{U_n}\exp\left(-{U_n \choose 2} \log\left(\frac{1}{t_n}\right)\right) \nonumber\\
&=& e^{-f_n U_n}, \nonumber
\end{eqnarray}
where~\(f_n\) is as defined in~(\ref{fn_def}). This proves the upper bound~(\ref{clq_upper}) in Proposition~\ref{clq_up_bound}.~\(\qed\)



In what follows, we estimate the probability~\(\mathbb{P}(B^c_L(S_q))\) to obtain the lower bounds in Theorem~\ref{clq_main_thm}.
We use the following Binomial estimate. Let~\(\{X_i\}_{1 \leq i \leq m}\) be independent Bernoulli random variables with~\[\mathbb{P}(X_i = 1) = p_i = 1-\mathbb{P}(X_i = 0).\] We have the following Lemma. 
\begin{Lemma}\label{bin_est_lem}Fix~\(0 < \epsilon  < \frac{1}{6}.\) If~\[T_m = \sum_{i=1}^{m} X_i,\] then
\begin{equation}\label{bin_est_f}
\mathbb{P}\left(|T_m - \mathbb{E}T_m| \geq \epsilon\mathbb{E}T_m\right) \leq \exp\left(-\frac{\epsilon^2 \mathbb{E}T_m}{4}\right)
\end{equation}
for all \(m \geq 1.\)
\end{Lemma}
For proof we refer to the Wikipedia link:\\ \(https://en.wikipedia.org/wiki/Chernoff\_bound\).
\subsection*{Small cliques estimate}
For integer~\(q \geq 1,\) we recall that~\(G(q,p)\) is the random graph with vertex set~\(S_q = \{1,2,\ldots,q\}.\) For integer~\(L \geq 2,\) let~\(B_L(S_q)\) denote the event that the random graph~\(G(q,\textbf{p})\) contains an open~\(L-\)clique; i.e., there are vertices~\(\{v_i\}_{1 \leq i \leq L}\) such that the edge between~\(v_i\) and~\(v_j\) is open for any~\(1 \leq i \neq j \leq L.\)
For~\(L \geq 2,\) we define
\begin{equation}\label{tlq}
t_L(q) = \sup_{\textbf{p} \in {\cal N}(a,q,p_q)}\mathbb{P}_{\textbf{p}}(B^c_L(S_q)).
\end{equation}
We first obtain a recursive relation involving~\(t_L(q).\)
\begin{Lemma}\label{cliq_rec_lem} Fix~\(0 < \epsilon < \frac{1}{6}\) and integer~\(L \geq 1.\) For integer~\(q \geq 1,\) suppose that~\(\textbf{p} \in {\cal N}(a,q,p_q)\)~(see~(\ref{beta_12def})) and let
\begin{equation}\label{q1_def}
q_1 = [(p_q-\delta)(q-1)]
\end{equation}
be the largest integer less than or equal to~\((p_q-\delta)(q-1).\) Let~\(\delta \in\{p_q\epsilon,(1-p_q)\epsilon\}.\) For all integers~\(q\) such that~\(q_1 \geq n^{a},\) we have that
\begin{eqnarray}
t_L(q) \leq  qt_{L-1}(q_1) + \exp\left(-\frac{\epsilon \delta}{10}q^2\right).\label{main_rec}
\end{eqnarray}
\end{Lemma}
\emph{Proof of Lemma~\ref{cliq_rec_lem}}: For simplicity, we write~\(p = p_q.\) We first prove that~(\ref{main_rec}) is satisfied with~\(\delta= p\epsilon.\) Let~\(N_e\) be the number of open edges in the random graph~\(G(q,\textbf{p}).\) Using~(\ref{beta_12def}), we have that~\(\mathbb{E}N_e \geq p{q \choose 2}.\) Fixing~\(0 < \epsilon <\frac{1}{6}\) and applying the binomial estimate~(\ref{bin_est_f}) with~\(T_m = N_e,\) we have that
\begin{equation}\label{n_e_est2}
\mathbb{P}_{\textbf{p}}\left(N_e \geq  p(1-\epsilon) {q \choose 2}\right) \geq  1- \exp\left(-\frac{\epsilon^2  p}{4} {q \choose 2}\right) = 1- \exp\left(-\frac{\epsilon\delta}{4} {q \choose 2}\right)
\end{equation}
for all~\(q \geq 2.\) The final term is obtained using~\(\delta = p\epsilon.\) Using~\(\frac{1}{4}{q \choose 2} \geq \frac{q^2}{10}\) for all~\(q \geq 5\) for the final term above we have
\begin{equation}\label{n_e_est}
\mathbb{P}_{\textbf{p}}\left(N_e \geq (p-\delta) {q \choose 2}\right) \geq 1-\exp\left(-\frac{\epsilon \delta}{10}q^2\right).
\end{equation}


Using~(\ref{n_e_est}), we therefore have
\begin{eqnarray}
\mathbb{P}_{\textbf{p}}(B^c_L(S_q)) = I_1 + I_2,\label{bl_est1}
\end{eqnarray}
where
\begin{equation}
I_1 := \mathbb{P}_{\textbf{p}}\left(B^c_L(S_q) \bigcap \left\{N_e \geq (p-\delta) {q \choose 2}\right\}\right) \label{i1_def}
\end{equation}
and
\begin{eqnarray}
I_2 &=& \mathbb{P}_{\textbf{p}}\left(B^c_L(S_q) \bigcap \left\{N_e < (p-\delta) {q \choose 2}\right\}\right)  \nonumber\\
&\leq& \mathbb{P}_{\textbf{p}}\left(N_e < (p-\delta) {q \choose 2}\right) \nonumber\\
&\leq& \exp\left(-\frac{\epsilon \delta}{10}q^2\right).\label{w4}
\end{eqnarray}



We estimate~\(I_1\) as follows. Suppose that the event~\(N_e \geq (p-\delta) {q \choose 2}\) occurs. If~\(d(v)\) denotes the degree of vertex~\(v \in \{1,2,\ldots,q\}\) in the random graph~\(G(q,\textbf{p}),\) we then have~\[\sum_{1 \leq v \leq q} d(v) = 2N_e \geq (p-\delta)q(q-1).\] In particular, there exists a vertex~\(w\) such that
\begin{equation}\label{d_w_est}
d(w) \geq (p-\delta)(q-1) \geq q_1.
\end{equation}
Here~\(q_1 = [(p-\delta)(q-1)]\) is as defined in the statement of the Lemma. This implies that
\begin{eqnarray}
&&\mathbb{P}_{\textbf{p}}\left(B^c_L(S_q) \bigcap \left\{N_e \geq (p-\delta){q \choose 2}\right\}\right) \nonumber\\
&&\;\;\;\;\;\;\leq \mathbb{P}_{\textbf{p}}\left(B^c_L(S_q) \bigcap \left(\bigcup_{1 \leq z \leq q}\left\{d(z) \geq q_1\right\}\right)\right) \nonumber\\
&&\;\;\;\;\;\;\leq \sum_{1 \leq z \leq q} \mathbb{P}_{\textbf{p}}\left(B^c_L(S_q) \bigcap \left\{d(z) \geq q_1\right\}\right).\label{q_est1}
\end{eqnarray}


Fixing~\(1 \leq z \leq q,\) we evaluate each term in~(\ref{q_est1}) separately. Letting~\(N(z) = N(z,G(q,\textbf{p}))\) be the set of neighbours of~\(z\) in the random graph~\(G(q,p),\) we have
\begin{eqnarray}
&&\mathbb{P}_{\textbf{p}}\left(B^c_L(S_q) \bigcap \left\{d(z) \geq q_1\right\}\right) \nonumber\\
&&\;\;\;\;\;\;\;\;\;\;\;= \sum_{S : \#S \geq q_1,\;z \notin S} \mathbb{P}_{\textbf{p}}\left(B^c_L(S_q) \bigcap \left\{N(z) = S \right\}\right). \label{b_l1}
\end{eqnarray}
Suppose now that the event~\(B^c_L(S_q)  \bigcap \{N(z) = S\}\) occurs for some fixed set~\(S\) with~\(\#S \geq q_1.\) We recall that since~\(B_L^{c}(S_q)\) occurs, there is no~\(L-\)clique in the random graph~\(G(q,\textbf{p})\) with vertex set~\(S_q = \{1,2,\ldots,q\}.\) This means that there is no~\((L-1)-\)clique in the random induced subgraph of~\(G(q,\textbf{p})\) formed by the vertices of~\(S;\) i.e., the event~\(B^{c}_{L-1}(S)\) occurs. Therefore we have
\begin{eqnarray}
\mathbb{P}_{\textbf{p}}\left(B^c_L(S_q) \bigcap \{N(z) = S\}\right) &\leq& \mathbb{P}_{\textbf{p}}\left(\{N(z) = S\} \cap B^c_{L-1}(S)\right) \nonumber\\
&=&\mathbb{P}_{\textbf{p}}\left(N(z) = S\right)\mathbb{P}_{\textbf{p}}\left(B^c_{L-1}(S)\right). \label{b_l2}
\end{eqnarray}
The equality~(\ref{b_l2}) true as follows. The event that~\(\{N(z) = S\}\) depends only on the state of edges containing~\(z\) as an endvertex. On the other hand, the event~\(B^c_{L-1}(S)\) depends only on the state of edges having both their endvertices in~\(S.\) Since the set~\(S\) does not contain the vertex~\(z\)~(see (\ref{b_l1})), we have that the events~\(\{N(z) = S\}\) and~\(B^c_{L-1}(S)\) are independent. This proves~(\ref{b_l2}).

We obtain the desired recursion using~(\ref{b_l2}) as follows. We recall that the set~\(S\) contains at least~\(q_1\) vertices (see~(\ref{b_l1})).
Therefore, setting~\(T\) to be the set of the~\(q_1\) least indices in~\(S,\) we have that if~\(B^c_{L-1}(S)\) occurs, then~\(B^c_{L-1}(T)\) occurs; i.e., there is no~\((L-1)-\)clique in the random induced subgraph formed by the vertices of~\(T.\) From~(\ref{b_l2}), we therefore have that
\begin{eqnarray}
\mathbb{P}_{\textbf{p}}\left(B^c_L(S_q) \bigcap \{N(z) = S\}\right) &\leq& \mathbb{P}_{\textbf{p}}\left(N(z) = S\right)\mathbb{P}_{\textbf{p}}\left(B^c_{L-1}(T)\right) \nonumber\\
&\leq&\mathbb{P}_{\textbf{p}}\left(N(z) = S\right)t_{L-1}(q_1). \label{b_l3}
\end{eqnarray}
The final inequality is true as follows. Let~\(\textbf{p}_T\) be the vector formed by the probabilities~\(\{p(i,j)\}_{i,j \in T}.\) From~(\ref{beta_12def}), we then have
\begin{equation}\label{beta_R}
\inf_{i \in T}\inf_{S} \frac{1}{\#S}\sum_{j \in S} p(i,j) \geq p
\end{equation}
for all~\(n \geq N.\) As in~(\ref{beta_12def}), the infimum is taken over all sets~\(S \subset T\) such that~\(\#S \geq n^{a}\) and~\(i \notin S.\) This proves that~\(\textbf{p}_T \in {\cal N}(a,q_1,p)\) and so~(\ref{b_l3}) is true.


Substituting~(\ref{b_l3}) into~(\ref{b_l1}), we have
\begin{eqnarray}
&&\mathbb{P}_{\textbf{p}}\left(B^c_L(S_q) \bigcap \left\{d(z) \geq q_1\right\}\right) \nonumber\\
&&\;\;\;\;\;\;\;\;\leq \sum_{S : \#S \geq q_1,\;z \notin S} \mathbb{P}_{\textbf{p}}\left(N(z) = S\right)t_{L-1}(q_1)\nonumber\\
&&\;\;\;\;\;\;\;\;=\mathbb{P}_{\textbf{p}}\left(\left\{d(z) \geq q_1\right\}\right) t_{L-1}(q_1)\label{b_l5}\\
&&\;\;\;\;\;\;\;\;\leq t_{L-1}(q_1). \label{b_l6}
\end{eqnarray}
The equality~(\ref{b_l5}) is true since the events~\(\{N(z) = S\}\) are disjoint for distinct~\(S.\)
Substituting~(\ref{b_l6}) into~(\ref{q_est1}), we have
\begin{eqnarray}
\mathbb{P}_{\textbf{p}}\left(B^c_L(S_q) \bigcap \left\{N_e \geq (p-\delta){q \choose 2}\right\}\right)
\leq \sum_{1 \leq z \leq q} t_{L-1}(q_1) = q t_{L-1}(q_1). \label{q_est2}
\end{eqnarray}

Using estimates~(\ref{q_est2}) and~(\ref{w4}) in~(\ref{bl_est1}) gives
\begin{eqnarray}
\mathbb{P}_{\textbf{p}}(B^c_L(S_q)) \leq  qt_{L-1}(q_1) + \exp\left(-\frac{\epsilon \delta}{10}q^2\right),\label{rec}
\end{eqnarray}
for all~\(q\) such that~\(q_1 \geq n^{a}.\) Taking supremum over all~\(\textbf{p} \in {\cal N}(a,q,p)\) proves~(\ref{main_rec}) with~\(\delta = p\epsilon.\)



It remains to see that~(\ref{main_rec}) is satisfied with~\(\delta = \epsilon(1-p).\) We recall that~\(N_e\) denotes the number of open edges in the random graph~\(G(n,p)\) (see the first paragraph of this proof). Let~\(W_e = {n \choose 2} - N_e\) denote the number of closed edges. Fixing~\(0 < \epsilon <\frac{1}{6}\) and applying the binomial estimate~(\ref{bin_est_f}) with~\(T_m = W_e,\) we have that
\begin{eqnarray}
\mathbb{P}_{\textbf{p}}\left(W_e \leq (1-p)(1+\epsilon) {q \choose 2}\right) &\geq&  1- \exp\left(-\frac{\epsilon^2(1-p)}{4} {q \choose 2}\right)
\nonumber\\
&=& 1- \exp\left(-\frac{\epsilon \delta}{4} {q \choose 2}\right) \nonumber
\end{eqnarray}
for all~\(q \geq 2.\) The final estimate follows using~\(\delta = \epsilon(1-p).\) Since
\[\left\{W_e \leq (1-p)(1+\epsilon) {q \choose 2} \right\} = \left\{N_e \geq (p-\delta) {q \choose 2} \right\},\] we again obtain~(\ref{n_e_est2}). The rest of the proof is as above. \(\qed\)

We use the recursion in the above Lemma iteratively to estimate the probability~\(t_L(q)\) of the event that there is no open~\(L-\)clique in the random graph~\(G(q,\textbf{p}).\)

\begin{Lemma}\label{t_q_lem} For integer~\(i \geq 1,\) define
\begin{equation}\label{v_i_def}
v_i = v_i(q) = (p-\delta)^{i}q - \frac{1}{1-p+\delta}.
\end{equation}
For all~\(q \geq 1\) such that~\(v_{L}(q) \geq n^{a},\) we have
\begin{equation}\label{t_sum}
t_L(q) \leq e^{-A_1} +2e^{-A_2}
\end{equation}
where
\begin{equation}\label{a1_def}
A_1 =  -L\log{q} + \log\left(\frac{1}{1-p}\right) \frac{v_{L}^2}{4}.
\end{equation}
and
\begin{equation}\label{a2_def}
A_2 =  \frac{\epsilon\delta}{10} v^2_{L} - L\log{q}.
\end{equation}
\end{Lemma}
To prove the above Lemma, we have a couple of preliminary estimates.
Let~\(\{q_i\}_{0 \leq i \leq L}\) be integers defined recursively as follows. The term~\(q_0 = q\) and for~\(i \geq 1,\) let \[q_i = [(p-\delta)(q_{i-1}-1)].\] For a fixed~\(1 \leq i \leq L,\) we have the following estimates.\\ 
\((a1)\) We have
\begin{equation}\label{one_step}
(p-\delta)(q_{i-1} -1 ) -1 \leq q_i \leq (p-\delta)q_{i-1} \leq q_{i-1} \leq q.
\end{equation}
\((a2)\) For~\(\delta > 0,\) we have
\begin{equation}\label{i_step}
v_i = (p-\delta)^{i}q - \frac{1}{1-p+\delta} \leq q_i \leq (p-\delta)^{i}q.
\end{equation}
\emph{Proof of~\((a1)-(a2)\)}: The property~\((a1)\) is obtained using the property~\(x-1 \leq [x] \leq x\) for any~\(x >0.\) Applying the upper bound in~(\ref{one_step}) recursively, we get \[q_i \leq (p-\delta)^{i}q_0 = (p-\delta)^{i} q.\] This proves the upper bound in~(\ref{i_step}). For the lower bound we again proceed iteratively and obtain for~\(i \geq 2\) that
\begin{eqnarray}
q_i &\geq& (p-\delta)(q_{i-1} -1 ) -1  \nonumber\\
&=& (p-\delta)q_{i-1} - ((p-\delta) + 1) \nonumber\\
&\geq& (p-\delta)^2q_{i-2}  - ((p-\delta)^2 + (p-\delta) + 1) \nonumber\\
&&\ldots \nonumber\\
&\geq&(p-\delta)^{i}q_{0} - \sum_{j=0}^{i}(p-\delta)^{j}. \label{q_i_def}
\end{eqnarray}
Since~\(\delta >0,\) we have \[\sum_{j=0}^{i}(p-\delta)^{j} \leq \frac{1}{1-p+\delta}\] and so
\[q_i \geq (p-\delta)^{i} q_0 - \frac{1}{1-p+\delta} = (p-\delta)^{i} q - \frac{1}{1-p+\delta}.\] This proves~\((a2).\)~\(\qed\)

Using the properties~\((a1)-(a2),\) we prove Lemma~\ref{t_q_lem}.\\
\emph{Proof of Lemma~\ref{t_q_lem}}: Letting
\begin{equation}\label{r_q_def}
r(q) = \exp\left(-\frac{\epsilon \delta}{10}q^{2}\right),
\end{equation}
we apply the recursion~(\ref{main_rec}) successively to get
\begin{eqnarray}
t_L(q) &\leq&  qt_{L-1}(q_1) + r(q)\nonumber\\
&\leq&  q (q_1 t_{L-2}(q_2) + r(q_1)) + r(q)\nonumber\\
&=&  q q_1 t_{L-2}(q_2) + qr(q_1) + r(q)\nonumber\\
&\leq&  q^2 t_{L-2}(q_2) + qr(q_1) + r(q)\nonumber
\end{eqnarray}
for all~\(q\) such that~\(q_2 = q_2(q) \geq n^{a}.\) The final estimate follows since~\(q_1 \leq q\) (see~(\ref{one_step}) of property~\((a1)\)).
Proceeding iteratively, we obtain the following estimate for all~\(q\) such that~\(q_{L-2}(q) \geq n^{a}:\)
\begin{eqnarray}
t_L(q) \leq  J_1 + J_2,\label{mult_time}
\end{eqnarray}
where
\begin{equation}\label{j1_def}
J_1 := q^{L-2}t_{2}(q_{L-2})
\end{equation}
and
\begin{equation}\label{j2_def}
J_2 := \sum_{j=0}^{L-3}q^{j}r(q_j).
\end{equation}

Let~\(v_L = v_L(q)\) be as defined in~(\ref{v_i_def}). For all~\(q\) such that~\(v_L(q) \geq n^{a},\) we have the following bounds for the terms~\(J_1\) and~\(J_2.\)
\begin{equation}\label{j1_est2}
J_1 \leq e^{-A_1}
\end{equation}
and~
\begin{equation}\label{j2_est2}
J_2 \leq 2e^{-A_2},
\end{equation}
where~\(A_1\) and~\(A_2\) are as given in~(\ref{a1_def}) and~(\ref{a2_def}), respectively. This proves the Lemma.\\
\emph{Proof of~(\ref{j1_est2}) and~(\ref{j2_est2})}: Since
\begin{equation}\label{qj_l2}
q_j \geq q_{L} \geq v_{L}
\end{equation}
for all~\(1 \leq j \leq L-1\) (property~\((a2)\)), the estimate~(\ref{mult_time}) holds for all~\(q\) such that~\(v_{L} = v_{L}(q) \geq n^{a}.\)

We first evaluate~\(J_1.\) We have
\begin{equation}
J_1 \leq q^{L-2}t_{2}(q_{L-2}) \leq q^{L}t_2(q_{L-2}) \label{j1_est1}
\end{equation}
and
\begin{equation}\label{two_cliq}
t_2(q_{L-2}) = (1-p)^{q_{L-2} \choose 2} \leq (1-p)^{q_{L} \choose 2} \leq (1-p)^{v_{L} \choose 2}.
\end{equation}
The first equality in~(\ref{two_cliq}) is true since there is no open \(2-\)clique among a set of vertices if and only if all the edges between the vertices are closed. The second and third inequality follow from~(\ref{qj_l2}) and the fact that~\(1-p < 1.\) Substituting~(\ref{two_cliq}) into~(\ref{j1_est1}) we get the estimate~(\ref{j1_est2}) for the term~\(J_1.\)

For the second term~\(J_2,\) we argue as follows. The term~\(r(q) = \exp\left(-\frac{\epsilon \delta}{10}q^2\right)\) defined in~(\ref{r_q_def}) is decreasing in~\(q.\) For~\(1 \leq j \leq L-1,\) we have from~(\ref{qj_l2}) that~\(q_j \geq v_L\) and so~\(r(q_j) \leq r(v_L).\) Using this in~(\ref{j2_def}), we then have
\begin{equation}
J_2 \leq \left(\sum_{j=0}^{L-3}q^{j}\right) r(v_{L})  = \frac{q^{L-2}-1}{q-1} r(v_{L}) \leq 2q^{L-3} r(v_{L}) \leq 2q^{L}r(v_{L}).\label{sec_term_t}
\end{equation}
The first inequality follows from the fact that~\(\frac{q^{L-2}-1}{q-1} \leq 2q^{L-3}\) for all~\(q \geq 3.\) Using the expression for~\(r(q)\) in~(\ref{sec_term_t}), we obtain~(\ref{j2_est2}).~\(\qed\)


\setcounter{equation}{0}
\renewcommand\theequation{\thesection.\arabic{equation}}
\section{Proof of Theorem~\ref{clq_main_thm}}\label{pf_clq_main_thm}
The following two estimates are used in what follows. For~\(0 < x < 1,\) we have
\begin{equation}\label{log_upper}
-\log(1-x) = \sum_{k\geq 1}\frac{x^{k}}{k} < \sum_{k\geq 1}x^{k} < \frac{x}{1-x},
\end{equation}
and
\begin{equation}\label{log_lower}
-\log(1-x) = \sum_{k\geq 1}\frac{x^{k}}{k}  > x.
\end{equation}

\subsection*{Proof of~\((i)\)}
Here~\(\alpha_1 > 0\) and we use the estimates~(\ref{a1_def}) and~(\ref{a2_def}) of Lemma~\ref{t_q_lem} to prove the Theorem~\ref{clq_main_thm}. We first obtain a couple of additional estimates. Fix~\(\eta\) and~\(\gamma\) as in the statement of the Theorem. Also fix~\(\epsilon  >0\) small to be determined later and set~\(q_n = n,\)
\begin{equation}\label{ln_def1}
L_n = (1-\eta) \frac{\log{n}}{\log\left(\frac{1}{p_n}\right)}
\end{equation}
and
\begin{equation}\label{deln_def1}
\delta_n = \epsilon p_n.
\end{equation}
For a fixed~\(\epsilon > 0,\) we have the following estimates regarding~\(L_n\) and~\(\delta_n.\)\\
\((b1)\) We have that
\begin{equation}\label{pn_raw1}
n^{-\alpha_1 - \epsilon} \leq p_n \leq n^{-\alpha_1+\epsilon} \text{ and } \frac{1}{1-p_n} \leq \frac{1}{1-n^{-\alpha_1+\epsilon}} \leq 2
\end{equation}
and so
\begin{equation}\label{pn_raw12}
\alpha_2 = \limsup_n \frac{\log\left(\frac{1}{1-p_n}\right)}{\log{n}} = 0.
\end{equation}
\((b2)\) We have
\begin{equation}\label{ln_est1}
L_n \leq \frac{1-\eta}{\alpha_1-\epsilon} \leq \frac{1}{\alpha_1-\epsilon}
\end{equation}
for all~\(n\) large.\\
\((b3)\) There is a constant~\(N_0 = N_0(\eta,\epsilon) \geq 1\) such that
\begin{equation}\label{vl_est1}
v_{L_n}  \geq n^{\eta -2\epsilon}
\end{equation}
for all~\(n\) large. \\\\
\emph{Proof of~\((b1)-(b3)\)}: We prove~\((b1)\) first. We use the definition of~\(\alpha_1 > 0\) to get that
\begin{equation}\label{pn_est_case1}
\frac{1}{\alpha_1+\epsilon} \leq \frac{\log{n}}{\log\left(\frac{1}{p_n}\right)} \leq \frac{1}{\alpha_1-\epsilon}
\end{equation}
for all~\(n\) large. This proves~(\ref{pn_raw1}) and~(\ref{pn_raw12}) in property~\((b1).\) 


The inequality in~(\ref{ln_est1}) follows from the final estimate of~(\ref{pn_est_case1}) and the definition of~\(L_n\) in~(\ref{ln_def1}). This proves~\((b2).\) To prove property~\((b3),\) we argue as  follows. Setting~\(q = q_n = n\) and~\(L=L_n\) in the definition of~\(v_i\) in~(\ref{v_i_def}), we then have
\begin{eqnarray}
v_{L_n} &=& \exp\left(L_n\log(p_n-\delta_n) + \log{n}\right) - \frac{1}{1-p_n+\delta_n} \nonumber\\
&\geq& \exp\left(L_n\log(p_n-\delta_n) + \log{n}\right) - \frac{1}{1-p_n} \nonumber\\
&\geq& e^{A_3} - 2 \label{vln_def1}
\end{eqnarray}
where
\begin{eqnarray}
A_3 &=& L_n\log(p_n-\delta_n)  + \log{n} \nonumber\\
&=& L_n \log\left(p_n(1-\epsilon)\right) + \log{n} \label{a_3_def1_c1_2}\\
&=& L_n\log{p_n} + L_n\log(1-\epsilon) + \log{n} \nonumber\\
&=& \eta \log{n} + L_n\log(1-\epsilon) \label{a_3_def1}
\end{eqnarray}
The final estimate in~(\ref{vln_def1}) follows from the final estimate in~(\ref{pn_raw1}). The equality~(\ref{a_3_def1_c1_2}) above is obtained using~\(\delta_n = \epsilon p_n\) and the final equality~(\ref{a_3_def1}) follows from the definition of~\(L_n\) in~(\ref{ln_def1}).

Using~\(-\log(1-x) < \frac{x}{1-x}\) (see~(\ref{log_upper})) with~\(x = \epsilon,\) we have
\[L_n\log(1-\epsilon) \geq \frac{-\epsilon}{1-\epsilon}L_n \geq \frac{-\epsilon}{1-\epsilon}\frac{1}{\alpha_1-\epsilon}\] where the final estimate follows using~(\ref{ln_est1}) in property~\((b2)\) above. Substituting the above into~(\ref{a_3_def1}), we have
\begin{eqnarray}\label{a3_est1_c1}
A_3 \geq \eta\log{n} - \frac{\epsilon}{1-\epsilon}\frac{1}{\alpha_1-\epsilon} \geq (\eta-\epsilon) \log{n}
\end{eqnarray}
for all~\(n\) large.
Using~(\ref{a3_est1_c1}) in~(\ref{vln_def1}), we have
\begin{equation}\label{vn_est33}
v_{L_n} \geq n^{\eta-\epsilon} - 2 \geq n^{\eta - 2\epsilon}
\end{equation}
for all~\(n\) large.~\(\qed\)

We use properties~\((b1)-(b3)\) to prove~\((i)\) in Theorem~\ref{clq_main_thm}. \\
\emph{Proof of \((i)\)}: From property~\((b3)\) and the choices of~\(\eta\) and~\(\gamma > 0\) as in the statement of the Thoerem, we have that~\(v_{L_n} \geq n^{a}\) for all~\(n \geq N_1\) large. Here~\(N_1 = N_1(\eta,\gamma,a)\) does not depend on the choice of~\(\textbf{p}.\) Thus the estimates for~\(A_1\) and~\(A_2\) in Lemma~\ref{t_q_lem} are applicable.

Setting~\(q = q_n = n\) and~\(L = L_n\) and~\(\delta = \delta_n\) as in~(\ref{ln_def1}) and~(\ref{deln_def1}), respectively, in the expressions for~\(A_1\) and~\(A_2\) in~(\ref{a1_def}) and~(\ref{a2_def}), we have
\begin{equation}\label{a1_est1_c1}
A_1 = A_1(n) = -L_n \log{n} + \log\left(\frac{1}{1-p_n}\right)\frac{v_{L_n}^2}{4},
\end{equation}
and
\begin{equation}\label{a2_est2_c1}
A_2 = A_2(n)=  \frac{\epsilon\delta_n}{10} v^2_{L_n} - L_n\log{n}.
\end{equation}
The following estimates for~\(A_1\) and~\(A_2\) imply the lower bound~(\ref{clq_loweri}) for case~\((i)\) in Theorem~\ref{clq_main_thm}. Fix~\(\gamma, \eta >0\) as in the statement of the Theorem.\\
\((c1)\) There are positive constants~\(\epsilon = \epsilon(\eta,\gamma) > 0\) and~\(M_1 = M_1(\eta,\gamma,\alpha_1) \geq 1\) so that
\begin{equation}\label{a1_est11_c1}
A_1 \geq n^{2\eta - 2\gamma-\alpha_1} 
\end{equation}
for all~\(n \geq M_1.\)\\ 
\((c2)\) There are positive constants~\(\epsilon = \epsilon(\eta,\gamma) > 0\) and~\(M_2 = M_2(\eta,\gamma,\alpha_1) \geq 1\) so that
\begin{equation}\label{a2_est22_c1}
A_2 \geq  n^{2\eta - 2\gamma-\alpha_1}
\end{equation}
for all~\(n \geq M_2.\)\\
\emph{Proof of~\((c1)-(c2)\)}: We first prove~\((c1).\) Using the estimate~(\ref{ln_est1}) of property~\((b2),\) we have that the first term in~(\ref{a1_est1_c1}) is
\begin{equation}\label{a11_est1_c1}
-L_n \log{n} \geq -\frac{\log{n}}{\alpha_1 -\epsilon}
\end{equation}
and using estimate~(\ref{vl_est1}) of property~\((b3),\) we have that the second term is
\begin{eqnarray}
\log\left(\frac{1}{1-p_n}\right)\frac{v_{L_n}^2}{4} &\geq& \log\left(\frac{1}{1-p_n}\right)\frac{n^{2\eta -4\epsilon}}{4} \nonumber\\
&\geq& p_n\frac{n^{2\eta -4\epsilon}}{4} \label{eq12_a1_c1}\\
&\geq& \frac{n^{2\eta - 5\epsilon - \alpha_1}}{4} \label{eq13_a1_c1}
\end{eqnarray}
for all~\(n\) large. The inequality~(\ref{eq12_a1_c1}) follows by setting~\(x = p_n\) in the estimate~\(-\log(1-x) > x\) (see~(\ref{log_lower})). The final estimate~(\ref{eq13_a1_c1}) follows from the first estimate~(\ref{pn_raw1}) of property~\((b1).\) 

Using estimates~(\ref{eq13_a1_c1}) and~(\ref{a11_est1_c1}) in the expression for~\(A_1\) in~(\ref{a1_est1_c1}), we have
\begin{eqnarray}
A_1 &\geq& \frac{1}{1-p_n}\left(\frac{n^{2\eta - 5\epsilon - \alpha_1}}{4} - \frac{\log{n}}{\alpha_1-\epsilon}\right) \nonumber\\
&\geq& \frac{1}{1-p_n}\left(\frac{n^{2\eta - 5\epsilon - \alpha_1}}{5}\right) \nonumber\\
&\geq& \frac{1}{5}n^{2\eta - 5\epsilon - \alpha_1} \nonumber
\end{eqnarray}
for all~\(n\) large. The final estimate follows using~\(1-p_n < 1.\)
We now fix~\(\gamma\) as in the statement of the Theorem and choose~\(\epsilon  = \epsilon(\gamma,\eta)> 0\) small so that
\[\frac{1}{5}n^{2\eta - 5\epsilon - \alpha_1} \geq n^{2\eta - 2\gamma-\alpha_1}\] for all~\(n\) large. This proves~\((c1).\)

We prove~\((c2)\) as follows. Using the upper bound for~\(L_n\log{n}\) in~(\ref{a11_est1_c1}) and the lower bound for~\(v_{L_n}\) in property~\((b3),\) we have
\begin{eqnarray}
A_2 &\geq& \frac{\epsilon\delta_n}{10} \frac{n^{2\eta-4\epsilon}}{4} - \frac{\log{n}}{\alpha_1 -\epsilon}  \nonumber\\
&\geq& \frac{\epsilon^2}{40} n^{2\eta-5\epsilon - \alpha_1} - \frac{\log{n}}{\alpha_1 -\epsilon} \label{a2_est211_c1}
\end{eqnarray}
where the final estimate~(\ref{a2_est211_c1}) follows from the fact that~\(\delta_n = \epsilon p_n\) and the lower bound for~\(p_n\) in~(\ref{pn_raw1}) (see property~\((b1)\)). As before, we choose~\(\epsilon  = \epsilon(\eta,\gamma) >0\) small so that the final term in~(\ref{a2_est211_c1}) is at least~\(n^{2\eta-2\gamma-\alpha_1}\) for all~\(n\) large. This proves~(\ref{a2_est2_c1}).~\(\qed\)


\subsection*{Proof of~\((ii)\)}
Fix~\(\eta\) and~\(\gamma\) as in the statement of the Theorem. Fix~\(\epsilon  >0\) small to be determined later and let~\(M = M(\epsilon) \geq 2\) be large so that
\begin{equation}\label{m_def}
\frac{1+\epsilon}{1-\frac{1+\epsilon}{M}} < 1+2\epsilon.
\end{equation}
Set~\(q_n = n,\)
\begin{equation}\label{ln_def2}
L_n = (1-\eta) \frac{\log{n}}{\log\left(\frac{1}{p_n}\right)}
\end{equation}
and
\begin{equation}\label{deln_def2}
\delta_n = \left(\epsilon_1 p_n \ind\left(p_n < 1- \frac{1}{M}\right) + \epsilon(1-p_n)\ind\left(p_n \geq 1-\frac{1}{M}\right)\right).
\end{equation}
Here~\(\epsilon_1 = \epsilon_1(\epsilon) > 0\) is to be determined later. For a fixed~\(\epsilon > 0,\) we have the following estimates regarding~\(L_n\) and~\(\delta_n.\)\\
\((b1)\) We have that
\begin{equation}\label{pn_raw2}
n^{-\epsilon } \leq p_n \leq 1 \text{ and } \frac{1}{1-p_n} \leq n^{\epsilon}.
\end{equation}
and
\begin{equation}\label{deln_c2}
\delta_n \geq n^{-2\epsilon}.
\end{equation}
for all~\(n\) large.\\
\((b2)\) We have
\begin{equation}\label{ln_est_c2}
L_n \leq (1-\eta)\frac{\log{n}}{1-p_n} \leq \frac{\log{n}}{1-p_n} \leq n^{\epsilon} \log{n}
\end{equation}
for all~\(n\) large. If~\(\epsilon_1 >0\) is sufficiently small, then
\begin{equation}\label{rat_bound}
R_n := \frac{\log\left(\frac{1}{p_n-\delta_n}\right)}{\log\left(\frac{1}{p_n}\right)} < 1+2\epsilon
\end{equation}
for all~\(n\) large.\\
\((b3)\) There is a constant~\(N_0 = N_0(\eta,\epsilon) \geq 1\) such that
\begin{equation}\label{vl_est2}
v_{L_n}  \geq n^{\eta -4\epsilon}
\end{equation}
for all~\(n\) large.\\\\
\emph{Proof of~\((b1)-(b3)\)}: The property~\((b1)\) is true as follows. Since~\(\alpha_1 = 0,\) we have from~(\ref{alp1_def}) that~\(\log\left(\frac{1}{p_n}\right) \leq \epsilon \log{n} = \log\left(n^{\epsilon}\right)\) for all~\(n\) large. This proves the first inequality in~(\ref{pn_raw2}). Since~\(\alpha_2 = 0,\) we have from~(\ref{alp2_def}) that~\(\log\left(\frac{1}{1-p_n}\right) \leq \epsilon \log{n} = \log\left(n^{\epsilon}\right)\) for all~\(n\) large. This proves the second inequality of~(\ref{pn_raw2}).

To prove~(\ref{deln_c2}), we proceed as follows. If~\(\delta_n = \epsilon_1 p_n,\) then we have from the first inequality in~(\ref{pn_raw2}) that~\(\delta_n \geq \epsilon_1 n^{-\epsilon} \geq n^{-2\epsilon}\) for all~\(n\) large. If~\(\delta_n = \epsilon(1-p_n),\) then using the second inequality in~(\ref{pn_raw2}), we have~\(\delta_n \geq \epsilon n^{-\epsilon} \geq n^{-2\epsilon}\) for all~\(n\) large.

The first estimate~(\ref{ln_est_c2}) follows by using the lower bound~\(-\log(1-x) > x\) with~\(x = 1-p_n\) in the definition of~\(L_n\) in~(\ref{ln_def2}). The second estimate in~(\ref{ln_est_c2}) follows using~\(\eta  < 1.\) The final estimate follows from~(\ref{pn_raw2}). To prove~(\ref{rat_bound}), we consider two cases separately depending on whether~\(\delta_n = \epsilon p_n\) or~\(\delta_n = (1-\epsilon)p_n.\)
If~\(\delta_n = \epsilon_1 p_n,\) then~\(p_n  < 1 - \frac{1}{M}\) and so we have
\[R_n = 1 + \frac{\log\left(\frac{1}{1-\epsilon_1}\right)}{\log\left(\frac{1}{p_n}\right)} \leq 1 + \frac{\log\left(\frac{1}{1-\epsilon_1}\right)}{\log\left(\frac{M}{M-1}\right)} \leq 1+\epsilon\] if~\(\epsilon_1 = \epsilon_1(\epsilon) > 0\) is small.

If~\(\delta_n = \epsilon (1-p_n),\) then~\(p_n \geq 1 - \frac{1}{M}\) and~\((1+\epsilon)(1-p_n) \leq \frac{1+\epsilon}{M} < 1\) since~\(M \geq 2\) and~\(0 < \epsilon  <1.\) Therefore using the upper bound estimate~\(-\log(1-x) < \frac{x}{1-x}\) from~(\ref{log_upper}) with~\(x = (1+\epsilon)(1-p_n),\) we have
\[-\log(p_n - \delta_n) = -\log\left(1 - (1+\epsilon)(1-p_n)\right) \leq \frac{(1+\epsilon)(1-p_n)}{1-(1+\epsilon)(1-p_n)}.\] Similarly using the lower bound estimate~\(-\log(1-x) > x\) from~(\ref{log_lower}), we have
\[-\log{p_n} = -\log(1-(1-p_n)) > 1-p_n.\] Using the above two estimates, we have
\[R_n \leq \frac{1+\epsilon}{1-(1+\epsilon)(1-p_n)} \leq \frac{1+\epsilon}{1 - \frac{1+\epsilon}{M}} \leq 1+2\epsilon\] by our choice of~\(M\) from~(\ref{m_def}).  This proves~\((b2).\)

To prove property~\((b3),\) we argue as  follows. Setting~\(q_n = n, L = L_n\) (as in~(\ref{ln_def2})) in the definition of~\(v_i\) in~(\ref{v_i_def}) we have
\begin{eqnarray}
v_{L_n} &=& \exp\left(L_n\log(p_n-\delta_n) + \log{n}\right) - \frac{1}{1-p_n+\delta_n} \nonumber\\
&\geq& \exp\left(L_n\log(p_n-\delta_n) + \log{n}\right) - \frac{1}{1-p_n} \nonumber\\
&=& \frac{1}{1-p_n}(e^{A_3} - 1) \label{vln_def2}
\end{eqnarray}
where
\begin{eqnarray}
A_3 &=& L_n\log(p_n-\delta_n) -\log\left(\frac{1}{1-p_n}\right)+ \log{n} \nonumber\\
&\geq& L_n\log(p_n-\delta_n) + (1-\epsilon)\log{n} \label{a_3_def2_e1}\\
&\geq& (-(1-\eta)(1+2\epsilon) + (1-\epsilon))\log{n} \label{a_3_def2_e2}\\
&=&\left(\eta(1+2\epsilon) - 3\epsilon\right)\log{n} \label{a3_def_fin}\\
&\geq&\left(\eta - 3\epsilon\right)\log{n} \label{a3_def_fin2}
\end{eqnarray}
for all~\(n\) large. The estimate in~(\ref{a_3_def2_e1}) follows since~\(\alpha_2\) defined in~(\ref{alp1_def}) is zero and so \(\log\left(\frac{1}{1-p_n}\right) < \epsilon \log{n}\) for all~\(n\) large. The estimate in~(\ref{a_3_def2_e2}) follows from~(\ref{rat_bound}) and the definition of~\(L_n\) in~(\ref{ln_def2}).

Substituting~(\ref{a3_def_fin2}) into~(\ref{vln_def2}), we get
\begin{equation}
v_{L_n} \geq \frac{1}{1-p_n}\left(n^{\eta - 3\epsilon} - 1\right) \geq n^{\eta-3\epsilon} - 1 \geq n^{\eta-4\epsilon}
\end{equation}
for all~\(n\) large. The second inequality follows using~\(1-p_n < 1.\) This proves~\((b3).\)\(\qed\)


We use properties~\((b1)-(b3)\) to prove~\((ii)\) in Theorem~\ref{clq_main_thm}. \\
\emph{Proof of~\((ii)\)}: We set~\(q = q_n = n\) and~\(L = L_n\) and~\(\delta = \delta_n\) as in~(\ref{ln_def2}) and~(\ref{deln_def2}), respectively, in the expressions for~\(A_1\) and~\(A_2\) in~(\ref{a1_def}) and~(\ref{a2_def}). We then have
\begin{equation}\label{a1_est1_c2}
A_1 = A_1(n) = -L_n \log{n} + \log\left(\frac{1}{1-p_n}\right)\frac{v_{L_n}^2}{4},
\end{equation}
and
\begin{equation}\label{a2_est2_c2}
A_2 = A_2(n)=  \frac{\epsilon\delta_n}{10} v^2_{L_n} - L_n\log{n}.
\end{equation}
The following estimates for~\(A_1\) and~\(A_2\) imply the lower bound~(\ref{clq_lowerii}) for case~\((ii)\) in Theorem~\ref{clq_main_thm}. Fix~\(\gamma, \eta >0\) as in the statement of the Theorem.\\
\((c1)\) There are positive constants \(\epsilon = \epsilon(\gamma,\eta) > 0\) and~\(M_1 = M_1(\gamma,\eta,\alpha_1,\epsilon) \geq 1\) so that
\begin{equation}\label{a1_est11_c2}
A_1 \geq n^{2\eta - 2\gamma}
\end{equation}
for all~\(n \geq M_1.\)\\ 
\((c2)\) There are positive constants~\(\epsilon  = \epsilon(\gamma,\eta) > 0\) and \(M_2 = M_2(\gamma,\eta,\alpha_1,\epsilon) \geq 1\) so that
\begin{equation}\label{a2_est22_c2}
A_2 \geq  n^{2\eta - 2\gamma}
\end{equation}
for all~\(n \geq M_2.\)\\
\emph{Proof of~\((c1)-(c2)\)}: We first prove~\((c1).\) Using the estimate~(\ref{ln_est_c2}) of property~\((b2),\) we have that the first term in~(\ref{a1_est1_c2}) is
\begin{equation}\label{a11_est_c2}
-L_n \log{n} \geq -n^{\epsilon}(\log{n})^2
\end{equation}
 and using estimate~(\ref{vl_est2}) of property~\((b3),\) we have that the second term in~(\ref{a1_est1_c2}) is
\begin{eqnarray}
\log\left(\frac{1}{1-p_n}\right)\frac{v_{L_n}^2}{4} &\geq& \log\left(\frac{1}{1-p_n}\right)\frac{n^{2\eta -8\epsilon}}{4} \nonumber\\
&\geq& p_n\frac{n^{2\eta -8\epsilon}}{4} \label{eq12_a1_c2}\\
&\geq& \frac{n^{2\eta - 9\epsilon}}{4} \label{eq13_a1_c2}
\end{eqnarray}
for all~\(n\) large. The inequality~(\ref{eq12_a1_c2}) follows using~\(-\log(1-x) > x\) for~\(0 < x< 1\) (see~(\ref{log_lower})). The final estimate in~(\ref{eq13_a1_c2}) follows from the first estimate~(\ref{pn_raw2}) of property~\((b1).\) 

Using estimates~(\ref{eq13_a1_c2}) and~(\ref{a11_est_c2}) in~(\ref{a1_est1_c2}), we have
\begin{eqnarray}
A_1 &\geq& \frac{n^{2\eta - 9\epsilon}}{4} - n^{\epsilon}(\log{n})^2\nonumber\\
&\geq& n^{2\eta-2\gamma}\label{a1_large2_c2}
\end{eqnarray}
for all~\(n\) large provided~\(\epsilon = \epsilon(\eta,\gamma) > 0\) is small. This proves~\((c1).\)

We prove~\((c2)\) as follows. Using the upper bound for~\(L_n\) in property~\((b2)\) and the lower bound for~\(v_{L_n}\) in property~\((b3),\) we have
\begin{eqnarray}
A_2 &\geq& \frac{\epsilon\delta_n}{10} \frac{n^{2\eta-8\epsilon}}{4} - n^{\epsilon}(\log{n})^2   \nonumber\\
&=& \frac{\epsilon}{40} n^{2\eta-10\epsilon} - n^{\epsilon}(\log{n})^2 \label{a2_est211_c2}\\
&\geq& n^{2\eta-2\gamma} \nonumber
\end{eqnarray}
for all~\(n\) large, provided~\(\epsilon = \epsilon(\eta,\gamma) > 0\) is small. The estimate~(\ref{a2_est211_c2}) follows from the estimate for~\(\delta_n\) in~(\ref{deln_c2}). This gives the estimate~\((c2)\) for the term~\(A_2.\)\(\qed\)


\subsection*{Proof of~\((iii)\)}
Fix~\(\eta\) and~\(\gamma\) as in the statement of the Theorem. Fix~\(\epsilon  >0\) small to be determined later and let~\(\alpha_2 = \alpha_2 > 0\) be as defined in~(\ref{alp2_def}). Set~\(q_n = n,\)
\begin{equation}\label{ln_def3}
L_n = (1-\eta) \frac{\log{n}}{\log\left(\frac{1}{p_n}\right)}
\end{equation}
and
\begin{equation}\label{deln_def3}
\delta_n = \epsilon(1-p_n).
\end{equation}
For a fixed~\(\epsilon > 0,\) we have the following estimates regarding~\(L_n\) and~\(\delta_n.\)\\
\((b1)\) We have that~\(\alpha_1 = 0\) and
\begin{equation}\label{pn_raw3}
n^{-\alpha_2 - \epsilon } \leq 1-p_n \leq n^{-\alpha_2+\epsilon} \text{ and } \frac{1}{p_n} \leq 2
\end{equation}
and
\begin{equation}\label{deln_est3}
\delta_n \geq n^{-\alpha_2-2\epsilon}
\end{equation}
for all~\(n\) large.\\
\((b2)\) We have
\begin{equation}\label{ln_est_c3}
L_n \leq (1-\eta-\alpha_2)\frac{\log{n}}{1-p_n} \leq \frac{\log{n}}{1-p_n} \leq n^{\alpha_2+\epsilon} \log{n}
\end{equation}
and
\begin{equation}\label{rat_bound_c3}
R_n := \frac{\log\left(\frac{1}{p_n-\delta_n}\right)}{\log\left(\frac{1}{p_n}\right)} < 1+2\epsilon
\end{equation}
for all~\(n\) large.\\
\((b3)\) We have that
\begin{equation}\label{vl_est3}
v_{L_n}  \geq n^{\eta +\alpha_2 - 5\epsilon}
\end{equation}
for all~\(n\) large.\\\\
\emph{Proof of~\((b1)-(b3)\)}: The property~\((b1)\) is true as follows. From the definition of~\(\alpha_2 > 0\) in~(\ref{alp2_def}) we have that
\begin{equation}\label{pn_est_case3}
\alpha_2 - \epsilon \leq \frac{\log\left(\frac{1}{1-p_n}\right)}{\log{n}} \leq \alpha_2 + \epsilon
\end{equation}
for all~\(n\) large. This proves the first inequality in~(\ref{pn_raw3}). The second inequality follows from the first inequality since
\[\frac{1}{p_n} \leq \frac{1}{1-n^{-\alpha_2+\epsilon}} \leq 2\] for all~\(n\) large. This also proves that~\(\alpha_1\) defined in~(\ref{alp1_def}) is zero. This proves the estimate~(\ref{pn_raw3}) of property~\((b1).\) To prove~(\ref{deln_est3}), we use~(\ref{pn_raw3}) and obtain \[\delta_n = \epsilon(1-p_n) \geq \epsilon n^{-\alpha_2-\epsilon} \geq n^{-\alpha_2-2\epsilon}\] for all~\(n\) large. This proves~\((b1).\)

The first estimate~(\ref{ln_est_c3}) follows by using the lower bound~\(-\log(1-x) > x\) with~\(x = 1-p_n\) in the definition of~\(L_n\) in~(\ref{ln_def2}). The second estimate in~(\ref{ln_est_c2}) follows using~\(1-\eta -\alpha_2 <1.\) The final estimate follows from~(\ref{pn_raw3}). The proof of~(\ref{rat_bound_c3}) is analogous as the proof of~(\ref{rat_bound}) for the case~\(\delta_n = \epsilon(1-p_n).\) This proves~\((b2).\)


To prove property~\((b3),\) we argue as  follows. Setting~\(q_n = n, L = L_n\) (as in~(\ref{ln_def3})) in the definition of~\(v_i\) in~(\ref{v_i_def}) we have
\begin{eqnarray}
v_{L_n} &=& \exp\left(L_n\log(p_n-\delta_n) + \log{n}\right) - \frac{1}{1-p_n+\delta_n} \nonumber\\
&=& \frac{1}{1-p_n+\delta_n}(e^{A_3} - 1) \label{vln_def3}
\end{eqnarray}

where
\begin{eqnarray}
A_3 &=& L_n\log(p_n-\delta_n) -\log\left(\frac{1}{1-p_n+\delta_n}\right) + \log{n} \nonumber\\
&=& L_n\log(p_n-\delta_n) +\log\left((1+\epsilon)(1-p_n)\right) + \log{n}. \label{a3_temp_c3}
\end{eqnarray}
The final equality is true using~\(\delta_n = \epsilon(1-p_n).\) For the middle term, we use the lower bound for~\(1-p_n\) from~(\ref{pn_raw3}) to get
\begin{equation}
\log\left((1+\epsilon)(1-p_n)\right) \geq \log(1-p_n) \geq  - (\alpha_2+\epsilon)\log{n}. \label{lg_first}\\
\end{equation}

We evaluate the first term in~(\ref{a3_temp_c3}) as follows. Since~\(\alpha_2 < 1,\) we have using~(\ref{rat_bound_c3}) and the definition of~\(L_n\) in~(\ref{ln_def3}) that
\begin{equation}\label{temp1_c3}
L_n\log(p_n-\delta_n) \geq -(1-\eta)(1+2\epsilon)\log{n}.
\end{equation}
Substituting~(\ref{temp1_c3}) and~(\ref{lg_first}) into~(\ref{a3_temp_c3}), we have
\begin{eqnarray}
A_3 &\geq&  -(1-\eta -\alpha_2)(1+2\epsilon)\log{n} -(\alpha_2+\epsilon)\log{n}  + \log{n}\nonumber\\
&=&(\eta(1+2\epsilon) - 3\epsilon + 2\epsilon\alpha_2)\log{n} \nonumber\\
&\geq&\left(\eta - 3\epsilon\right)\log{n} \label{a3_def_fin3}
\end{eqnarray}
for all~\(n\) large. Substituting~(\ref{a3_def_fin3}) into~(\ref{vln_def3}), we get
\begin{eqnarray} \nonumber
v_{L_n} &\geq& \frac{1}{1-p_n+\delta_n}\left(n^{\eta - 3\epsilon} - 1\right) \nonumber\\
&\geq& \frac{n^{\eta-4\epsilon} }{1-p_n+\delta_n} \nonumber\\
&\geq& \frac{n^{\eta-4\epsilon} }{1-p_n} \nonumber\\
&\geq& n^{\eta+\alpha_2-5\epsilon}
\end{eqnarray}
for all~\(n\) large. The final inequality follows from the estimate~(\ref{pn_raw3}) in property~\((b1).\) This proves~\((b3)\) for the case~\(\alpha_2 < 1.\)~\(\qed\)



We use properties~\((b1)-(b3)\) to prove~\((iii)\) in Theorem~\ref{clq_main_thm}. \\
\emph{Proof of~\((iii)\)}: We set~\(q = q_n = n\) and~\(L = L_n\) and~\(\delta = \delta_n\) as in~(\ref{ln_def3}) and~(\ref{deln_def3}), respectively, in the expressions for~\(A_1\) and~\(A_2\) in~(\ref{a1_def}) and~(\ref{a2_def}). We then have
\begin{equation}\label{a1_est1_c3}
A_1 = A_1(n) = -L_n \log{n} + \log\left(\frac{1}{1-p_n}\right)\frac{v_{L_n}^2}{4},
\end{equation}
and
\begin{equation}\label{a2_est2_c3}
A_2 = A_2(n)=  \frac{\epsilon\delta_n}{10} v^2_{L_n} - L_n\log{n}.
\end{equation}
The following estimates for~\(A_1\) and~\(A_2\) imply the lower bound~(\ref{clq_loweriii}) in case~\((iii)\) of Theorem~\ref{clq_main_thm}. Fix~\(\gamma, \eta >0\) as in the statement of the Theorem.\\
\((c1)\) There are positive constants \(\epsilon = \epsilon(\gamma,\eta) > 0\) and~\(M_1 = M_1(\gamma,\eta,\alpha_1,\epsilon) \geq 1\) so that
\begin{equation}\label{a1_est11_c3}
A_1 \geq n^{2\eta - 2\gamma + 2\alpha_2}
\end{equation}
for all~\(n \geq M_1.\)\\ 
\((c2)\) There are positive constants~\(\epsilon  = \epsilon(\gamma,\eta) > 0\) and \(M_2 = M_2(\gamma,\eta,\alpha_1,\epsilon) \geq 1\) so that
\begin{equation}\label{a2_est22_c3}
A_2 \geq  n^{2\eta - 2\gamma + \alpha_2}
\end{equation}
for all~\(n \geq M_2.\)\\
\emph{Proof of~\((c1)-(c2)\)}: We first prove~\((c1).\) Using the estimate~(\ref{ln_est_c3}) of property~\((b2),\) we have that the first term in~(\ref{a1_est1_c3}) is
\begin{equation}\label{a11_est_c3}
-L_n \log{n} \geq -n^{\alpha_2 + \epsilon}(\log{n})^2
\end{equation}
and using estimate~(\ref{vl_est3}) of property~\((b3),\) we have that the second term in~(\ref{a1_est1_c3}) is
\begin{eqnarray}
\log\left(\frac{1}{1-p_n}\right)\frac{v_{L_n}^2}{4} &\geq& \log\left(\frac{1}{1-p_n}\right)\frac{n^{2\eta +2\alpha_2 -10\epsilon}}{4} \nonumber\\
&\geq& p_n\frac{n^{2\eta + 2\alpha_2-10\epsilon}}{4} \label{eq12_a1_c3}\\
&\geq& \frac{n^{2\eta +2\alpha_2 - 10\epsilon}}{8} \label{eq13_a1_c3}
\end{eqnarray}
for all~\(n\) large. The inequality~(\ref{eq12_a1_c3}) follows using~\(-\log(1-x) > x\) for~\(0 < x< 1\) (see~(\ref{log_lower})). The final estimate in~(\ref{eq13_a1_c3}) follows from the final estimate~(\ref{pn_raw3}) of property~\((b1).\) 

Using estimates~(\ref{eq13_a1_c3}) and~(\ref{a11_est_c3}) in~(\ref{a1_est1_c3}), we have
\begin{eqnarray}
A_1 &\geq& \frac{n^{2\eta +2\alpha_2 - 10\epsilon}}{8} - n^{\alpha_2+\epsilon}(\log{n})^2\nonumber\\
&\geq& n^{2\eta-2\gamma + 2\alpha_2}\label{a1_large2_c3}
\end{eqnarray}
for all~\(n\) large, provided~\(\epsilon = \epsilon(\eta,\gamma) > 0\) is small. This proves~\((c1).\)

We prove~\((c2)\) as follows. Using the upper bound for~\(L_n\) in property~\((b2)\) and the lower bound for~\(v_{L_n}\) in property~\((b3),\) we have
\begin{eqnarray}
A_2 &\geq& \frac{\epsilon\delta_n}{10} \frac{n^{2\eta+2\alpha_2-10\epsilon}}{4} - n^{\alpha_2+\epsilon}(\log{n})^2   \nonumber\\
&=& \frac{\epsilon}{40} n^{2\eta+\alpha_2-12\epsilon} - n^{\alpha_2 + \epsilon}(\log{n})^2 \label{a2_est211_c3}\\
&\geq& n^{2\eta-2\gamma + \alpha_2} \nonumber
\end{eqnarray}
for all~\(n\) large, provided~\(\epsilon = \epsilon(\eta,\gamma) > 0\) is small. The estimate~(\ref{a2_est211_c3}) follows from the estimate for~\(\delta_n\) in~(\ref{deln_est3}). This gives the estimate~\((c2)\) for the term~\(A_2.\)\(\qed\)

\setcounter{equation}{0}
\renewcommand\theequation{\thesection.\arabic{equation}}
\section{Proof of Proposition~\ref{clq_extrem_lem} and Theorem~\ref{clq_cor}} \label{pf_clq_cor}
\emph{Proof of Proposition~\ref{clq_extrem_lem}}: By definition of~\(\alpha_2\) in~(\ref{alp2_def}), we have \[(\alpha_2-\epsilon)\log{n} \leq \log\left(\frac{1}{1-p_n}\right) \leq (\alpha_2 +\epsilon) \log{n}\] so that
\begin{equation}\label{eq_pn_ext1}
\frac{1}{n^{\alpha_2+\epsilon}} \leq 1-p_n \leq \frac{1}{n^{\alpha_2-\epsilon}}
\end{equation}
for all~\(n\) large. Since~\(\alpha_2 > 2,\) we fix~\(\epsilon > 0\) small so that~\(\alpha_2 -\epsilon  > 2.\) If~\(N_e\) denote the number of open edges in the random graph~\(G(n,p_n),\) we then have
\begin{equation}\label{ne_est_alpha2}
\mathbb{P}\left(N_e \geq 1\right) \leq \mathbb{E} N_e  = p_n{n \choose 2} \leq \frac{n^2}{n^{\alpha_2-\epsilon}} \longrightarrow 0
\end{equation}
as~\(n \rightarrow \infty.\) But~\(\{N_e = 0\} = \{\omega(G(n,p_n)) = 1\}\) and so we obtain~(\ref{clq_extrem1}).

An analogous proof holds for the other case~\(\alpha_2 > 2\) by considering closed edges. 

If~\(1 < \alpha_2 < 2,\) we argue as follows. If~\(W_e\) denotes the number of closed edges, then using the Binomial estimate~(\ref{bin_est_f}), we have
\begin{equation}\label{eq_we_ext1}
\mathbb{P}\left(\left|W_e -\mathbb{E} W_e\right| \geq \epsilon \mathbb{E} W_e\right) \leq \exp\left(-\frac{\epsilon^2(1-p_n)}{4}{n \choose 2}\right).
\end{equation}
Using~(\ref{eq_pn_ext1}), we have
\begin{equation}\label{eq_we_upp}
\mathbb{E}W_e  = (1-p_n){n \choose 2} \leq \frac{1}{n^{\alpha_2-\epsilon}}\frac{n^2}{2} \leq \frac{1}{2}n^{2-\alpha_2+\epsilon}
\end{equation}
and
\begin{equation}\label{eq_we_low}
\mathbb{E}W_e  = (1-p_n){n \choose 2} \geq \frac{1}{n^{\alpha_2+\epsilon}}\frac{n^2}{4} \geq \frac{1}{4}n^{2-\alpha_2-\epsilon}
\end{equation}
for all~\(n\) large. The first inequality in~(\ref{eq_we_low}) is obtained using~(\ref{eq_pn_ext1}) and~\({n \choose 2} \geq \frac{n^2}{4}\) for all~\(n\) large. We choose~\(\epsilon >0\) small so that \[0 < 2-\alpha_2 - 2\epsilon < 2-\alpha_2+2\epsilon < 1.\]

We then have from~(\ref{eq_we_ext1}),~(\ref{eq_we_low}) and~(\ref{eq_we_upp}) that
\begin{eqnarray}\label{eq_we_ext2}
\mathbb{P}\left(W_e \geq (1+\epsilon) \frac{1}{2}n^{2-\alpha_2+\epsilon}\right) &\leq& \mathbb{P}\left(W_e \geq (1+\epsilon) \mathbb{E} W_e\right) \nonumber\\
&\leq& \exp\left(-\frac{\epsilon^2}{4}\frac{1}{4}n^{2-\alpha_2-\epsilon}\right) \nonumber\\
&\leq& \exp\left(-n^{2-\alpha_2-2\epsilon}\right) \nonumber
\end{eqnarray}
for all~\(n\) large. Suppose now that the event~\(W_e \leq (1+\epsilon) \frac{1}{2}n^{2-\alpha_2+\epsilon}\) occurs and let~\({\cal S}_e\) be the set of all vertices belonging to the closed edges in the random graph~\(G(n,p_n).\) The induced subgraph~\(G_S\) with vertex set~\(\{1,2,\ldots,n\} \setminus {\cal S}_e\) contains at least~\(n-(1+\epsilon)n^{2-\alpha_2+\epsilon}\) vertices and every edge in~\(G_S\) is open. In other words, the graph~\(G_S\) is an open clique containing at least~\[n-(1+\epsilon)n^{2-\alpha_2+\epsilon} \geq n-n^{2-\alpha_2+2\epsilon}\] vertices, for all~\(n\) large.

We now prove the upper bound~(\ref{clq_upper}). For integer~\(q \geq 1,\) let~\(G(q,p)\) be the random graph with vertex set~\(S_q = \{1,2,\ldots,q\}.\) For integer~\(L \geq 2,\) let~\(B_L(S_q)\) denote the event that the random graph~\(G(q,p)\) contains an open~\(L-\)clique; i.e., there are vertices~\(\{v_i\}_{1 \leq i \leq L}\) such that the edge between~\(v_i\) and~\(v_j\) is open for any~\(1 \leq i \neq j \leq L.\)
We have
\begin{equation}\label{blq_upper}
\mathbb{P}(B_L(S_q)) \leq {q \choose L}p^{L \choose 2} \leq q^{L}p^{L \choose 2} = e^{-LA_0}
\end{equation}
where
\begin{equation}\label{a0_def}
A_0 = A_0(q,p,L) = \left(\frac{L-1}{2}\right)\log\left(\frac{1}{p}\right) - \log{q}.
\end{equation}
We now set~\(q = n, p = p_n\) and let~\(f_n \rightarrow \infty\) be any sequence as in the statement of the Theorem.
Setting~\(L =U_n\) as defined in~(\ref{un_def}), we then have \(A_0 = A_0(n) = f_n.\) This proves the upper bound~(\ref{clq_upper}) in Lemma~\ref{clq_extrem_lem}.~\(\qed\)

\emph{Proof of~\((i)\)}: Here~\(\alpha_1\) defined in~(\ref{alp1_def}) equals~\(\theta_1\) and~\(\alpha_2\) as defined in~(\ref{alp2_def}) equals zero. The lower bound follows from~(\ref{clq_loweri}), case~\((i)\) of Theorem~\ref{clq_main_thm}. For the upper bound, we fix~\(\xi > 0\) and set~\(f_n = \xi \log{n}\) so that~\(U_n\) as defined in~(\ref{un_def}) equals~\(\frac{2+\xi}{\theta_1} + 1.\) The upper bound then follows from~(\ref{clq_upper}).~\(\qed\)

\emph{Proof of~\((ii)\)}: Here~\(\alpha_1\) and~\(\alpha_2\) defined in~(\ref{alp1_def}) and~(\ref{alp2_def}), respectively, both equal zero.  Fixing~\(\eta,\gamma\) as in the statement of~\((ii),\) the lower bound follows from~(\ref{clq_lowerii}), case~\((ii)\) of Theorem~\ref{clq_main_thm}.

For the upper bound, we fix~\(0 < \xi_1 < \xi < 1\) and set~\(f_n = \xi_1 \log{n}.\) The term~\(U_n\) defined in~(\ref{un_def}) equals
\[U_n = \frac{(2+\xi_1)\log{n}}{\log\left(\frac{1}{p}\right)} + 1 \leq \frac{(2+\xi)\log{n}}{\log\left(\frac{1}{p}\right)}\] for all~\(n \geq N_1.\) Here~\(N_1 = N_1(\xi,\xi_1,p) \geq 1\) is a constant. Using~(\ref{clq_upper}) of Theorem~\ref{clq_main_thm}, we have
\begin{equation}\label{clq_upper_22}
\mathbb{P}\left(\omega(G(n,p)) \leq \frac{(2+\xi)\log{n}}{\log\left(\frac{1}{p}\right)}\right) \geq 1-\exp\left(-\frac{\xi_1(2+\xi_1)}{\log\left(\frac{1}{p}\right)} (\log{n})^2\right)
\end{equation}
for all~\(n \geq N_1.\) Choosing~\(\xi_1\) sufficiently close to~\(\xi\) so that~\(\xi_1(2+\xi_1) >\xi(1+\xi),\) we obtain the upper bound in~\((ii).\)~\(\qed\)

\emph{Proof of~\((iii)\)}: Here~\(\alpha_1\) defined in~(\ref{alp1_def}) equals zero and~\(\alpha_2\) as defined in~(\ref{alp2_def}) equals~\(\theta_2.\) Let~\(\eta,\gamma> 0\) be as in the statement of the Theorem and fix~\(\gamma_0,\eta_0 > 0\) such that~\(\gamma_0 -\frac{\theta_2}{2}< \gamma -\frac{\theta_2}{2} <  \eta_0 < \eta < 1-\theta_2\) and~\(\eta_0-\gamma_0 > \eta-\gamma.\) Let~\(\epsilon >0\) be small to be determined later. Applying the lower bound~(\ref{clq_loweriii}), case~\((iii)\) in Theorem~\ref{clq_main_thm} with~\(\eta_1\) and~\(\gamma_1\) we have
\begin{eqnarray}
&&\mathbb{P}\left(\omega(G(n,p_n)) \geq (1-\theta_2-\eta_0)\frac{\log{n}}{\log\left(\frac{1}{p_n}\right)}\right) \nonumber\\
&&\;\;\;\;\;\geq 1 - 3\exp\left(-n^{2\eta_0 - 2\gamma_0+\theta_2}\right) \nonumber\\
&&\;\;\;\;\;\geq 1 - 3\exp\left(-n^{2\eta - 2\gamma+\theta_2}\right), \label{lower_case3}
\end{eqnarray}
where the final estimate follows from the choices of~\(\eta_0\) and~\(\gamma_0.\)
We have
\begin{equation}\label{log_pn_bound}
\log\left(\frac{1}{p_n}\right) = -\log(1-(1-p_n)) < \frac{1-p_n}{1-(1-p_n)} \leq \frac{1-p_n}{1-\epsilon} = \frac{1}{n^{\theta_2}(1-\epsilon)}
\end{equation}
for all~\(n \geq N_1.\) Here~\(N_1 = N_1(\epsilon) \geq 1\) is a constant. The first inequality in~(\ref{log_pn_bound}) follows from~(\ref{log_upper}) and the second inequality follows from the fact that~\(1-p_n < \epsilon\) for all~\(n \geq N_1\) large. From~(\ref{log_pn_bound}), we therefore have
\begin{eqnarray}
(1-\theta_2- \eta_0)\frac{\log{n}}{\log\left(\frac{1}{p_n}\right)} &\geq&  (1-\theta_2-\eta_0)(1-\epsilon)n^{\theta_2}\log{n}\nonumber\\
&\geq& (1-\theta_2-\eta)n^{\theta_2}\log{n},\label{term_bound1}
\end{eqnarray}
provided~\(\epsilon = \epsilon(\eta_1,\eta,\theta_2) > 0\) is small. Fixing such an~\(\epsilon\) and substituting the estimate~(\ref{term_bound1}) into~(\ref{lower_case3}), we obtain the lower bound in~(\ref{clq_bd_case3}).

For the upper bound, we fix~\(\xi >0\) and set~\(f_n = \xi \log{n}.\) The term~\(U_n\) defined in~(\ref{un_def}) is then
\begin{equation}\label{un_case3}
U_n = \frac{(2+\xi)\log{n}}{\log\left(\frac{1}{p_n}\right)}
\end{equation}
and using~(\ref{log_lower}), we have
\begin{equation}\label{log_pn_bound2}
\log\left(\frac{1}{p_n}\right) = -\log(1-(1-p_n))  > 1-p_n = \frac{1}{n^{\theta_2}}.
\end{equation}
Using the bounds~(\ref{log_pn_bound}) and~(\ref{log_pn_bound2}) in~(\ref{un_case3}), we have
\begin{equation}\label{un_case3_bound}
(2+\xi)(1-\epsilon)n^{\theta_2}\log{n} \leq U_n \leq (2+\xi)n^{\theta_2}\log{n}
\end{equation}
for all~\(n \geq N_1.\) Here~\(N_1 = N_1(\epsilon) \geq 1\) is the constant in~(\ref{log_pn_bound}).

Using the above bounds in the upper bound~(\ref{clq_upper}) of Theorem~\ref{clq_main_thm}, we have
\begin{equation}\label{clq_upper_cs3}
\mathbb{P}\left(\omega(G(n,p_n)) \leq (2+\xi)n^{\theta_2}\log{n}\right) \geq 1-\exp\left(-\xi(2+\xi)(1-\epsilon)n^{\theta_2}(\log{n})^2\right)
\end{equation}
for all~\(n \geq N_1.\) Choosing~\(\epsilon > 0\) small so that~\((2+\xi)(1-\epsilon) > 1+\xi,\) we obtain the upper bound in~(\ref{clq_bd_case3}).~\(\qed\)


\setcounter{equation}{0}
\renewcommand\theequation{\thesection.\arabic{equation}}
\section{Proof of Theorem~\ref{chr_cor}} \label{pf_chr_cor}
For a graph~\(G = (V,E)\) on~\(n\) vertices, let \(\alpha(G)\) be the independence number of the graph~\(G\) defined as as follows. For integer~\(0 \leq h \leq n,\) we say that~\(\alpha(G) = h\) if and only if the following two conditions are satisfied:\\ \((a)\) There is a set of~\(h\) vertices, none of which have an edge between them.\\ \((b)\) If~\(h+1 \leq n,\) then every set of~\(h+1\) vertices have an edge between them.

As in Section~\ref{intro}, let~\(\omega(G)\) denote the clique number of~\(G.\) Let~\(\overline{G} = (\overline{V},\overline{E})\) denote the compliment of the graph~\(G\) defined as follows. The vertex set~\(\overline{V} = V\) and an edge~\(e \in \overline{E}\) if and only if~\(e \notin E.\) The following three properties are used to prove Theorem~\ref{chr_cor}.\\
\((d1)\) We have
\begin{equation}\label{alpha_omega}
\alpha(G) = \omega(\overline{G}).
\end{equation}
\((d2)\) We have
\begin{equation}\label{chi_lower}
\chi(G) \geq \frac{n}{\alpha(G)} = \frac{n}{\omega(\overline{G})}.
\end{equation}
\((d3)\) Suppose for some integer~\(1 \leq m \leq n,\) every set of~\(m\) vertices in the complement graph~\(\overline{G}\) contains a clique of size~\(L.\) We then have
\begin{equation}\label{chi_upper}
\chi(G) \leq \frac{n-m}{L} + m +1 \leq \frac{n}{L} + 2m.
\end{equation}



The lower bounds in~Theorem~\ref{chr_cor}, follow from the respective upper\\bounds~(\ref{clq_bd_case1}),~(\ref{clq_bd_case2}) and~(\ref{clq_bd_case3}) on the clique number~\(\omega(G(n,1-r_n))\) of~Theorem~\ref{clq_cor} and property~\((d2)\) above. This is because, the random graph~\(\overline{G}(n,r_n)\) has the same distribution as the random graph~\(G(n,1-r_n).\)

For the upper bounds, we consider each case separately. \\\\
\emph{Proof of~\((i)\)}: Here~\(r_n = \frac{1}{n^{\theta_2}}\) for some~\(\theta_2 > 0.\) To estimate the chromatic number using property~\((d3),\) we identify cliques in subsets of the random graph~\(G(n,1-r_n).\) Fix~\(\beta > 0\) to be determined later and set~\(m = n^{1-\beta}\) and apply Theorem~\ref{clq_cor}, case~\((iii)\) for the random graph~\(G(m,p_m),\) where
\begin{equation}\label{pm_def}
p_m = 1-\frac{1}{m^{\theta_{22}}} = 1-r_n,
\end{equation}
where~\(\theta_{22} = \frac{\theta_2}{1-\beta}.\) We then have~\(\alpha_1 = 0\) and~\(\alpha_2 = \theta_{22},\) where~\(\alpha_1\) and~\(\alpha_2\) are as defined in~(\ref{alp1_def}) and~(\ref{alp2_def}), respectively. Let~\(\eta,\gamma >0\) be such that
\begin{equation}\label{chr_condi}
\frac{1-\theta_{22}}{2}  + \gamma < \eta <1-\theta_{22}
\end{equation}

From the proof of lower bound of~(\ref{clq_bd_case3}), there is a positive integer~\(N_3= N_3(\eta,\gamma,\xi) \geq 1\) so that
\begin{eqnarray}
\mathbb{P}\left(\omega(G(m,p_m)) \leq L\right) \leq 3\exp\left(-m^{2\eta - 2\gamma+\theta_{22}}\right) \label{clq_bd_case3m_c1}
\end{eqnarray}
for all~\(m\) large, where
\begin{equation}\label{ldef}
L = (1-\eta-\theta_{22})m^{\theta_{22}} \log{m} = (1-\eta-\theta_{22})(1-\beta)n^{\theta_{2}} \log{n}.
\end{equation}
The final estimate above follows using~\(m = n^{1-\beta}.\)

Let~\({\cal S}_m\) be the set of subsets of size~\(m\) in~\(\{1,2,\ldots,n\}\) and for a set~\(S \in {\cal S}_m,\) let~\(F_n(S)\) denote the event that the random induced subgraph of~\(G(n,1-r_n)\) with vertex set~\(S\) contains an open~\(L-\)clique. From~(\ref{clq_bd_case3m_c1}) we have that
\[\mathbb{P}(F^c_n(S)) \leq 3\exp\left(-m^{2\eta - 2\gamma-\theta_{22}}\right) \] for all~\(n\) large. Let~\(F_n = \bigcap_{S \in {\cal S}_m}F_n(S)\) denote the event that every set of~\(m\) vertices in the random graph~\(G(n,1-r_n)\) contains an~\(L-\)clique. Since there are~\({n \choose m}\) sets in~\({\cal S}_m,\) we have
\begin{equation}\label{fn_est}
\mathbb{P}(F_n^c) \leq {n \choose m}3\exp\left(-m^{2\eta - 2\gamma+\theta_{22}}\right) \leq n^{m} 3\exp\left(-m^{2\eta - 2\gamma+\theta_{22}}\right) = 3e^{-B},
\end{equation}
where
\begin{equation}\label{best_c1}
B = m^{2\eta - 2\gamma-\theta_{22}} - m\log{n} = m^{2\eta - 2\gamma-\theta_{22}} - \frac{m}{1-\beta}\log{m}.
\end{equation}
The final estimate follows since~\(m = n^{1-\beta}.\) From the choices of~\(\eta\) and~\(\gamma\) in~(\ref{chr_condi}), we have that~\(2\eta - 2\gamma+\theta_{22} > 1\) and so
\begin{equation}\label{besti}
B \geq \frac{1}{2}m^{2\eta-2\gamma+\theta_{22}} = \frac{1}{2}n^{(1-\beta)(2\eta-2\gamma)+\theta_2}
\end{equation}
for all~\(n \geq N_2.\) Here~\(N_2 = N_2(\eta,\gamma,\beta) \geq 1\) is a constant and the final equality follows from the definition of~\(m = n^{1-\beta}\) and~\(\theta_{22} = \frac{\theta_2}{1-\beta}\) above.

If the event~\(F_n\) occurs, then using property~\((d3),\) we have that
\begin{equation}\label{chi_uper_casei_est}
\chi(G(n,r_n)) \leq \frac{n}{L} + 2m  = \frac{1}{(1-\eta-\theta_{22})(1-\beta)}\frac{n^{1-\theta_2}}{\log{n}} + 2n^{1-\beta}
\end{equation}
Fixing~\(\beta > \theta_2\) and~\(\xi > 0,\) we have that the final expression in~(\ref{chi_uper_casei_est}) is at most~\[\left(\frac{1+0.5\xi}{(1-\eta-\theta_{22})(1-\beta)}\right)\frac{n^{1-\theta_2}}{\log{n}}\] for all~\(n\)  large.
Summarizing, we have from~(\ref{best_c1}),~(\ref{besti}) and~(\ref{fn_est}) that
\begin{eqnarray}
&&\mathbb{P}\left(\chi(G(n,r_n)) \leq \left(\frac{1+0.5\xi}{(1-\eta-\theta_{22})(1-\beta)}\right)\frac{n^{1-\theta_2}}{\log{n}}\right) \nonumber\\
&&\;\;\;\;\;\geq 1 - 3\exp\left(-\frac{1}{2}n^{(1-\beta)(2\eta-2\gamma) + \theta_2}\right).\label{chi_est_c1}
\end{eqnarray}

We have the following property.\\
\((f1)\) Let \[{\cal T} = \{(\eta,\gamma,\beta) : \beta > \theta_2 \text{ and~(\ref{chr_condi}) is satisfied}\}.\]
Let~\(\xi,\zeta > 0\) be as in the statement of the Theorem. There exists~\((\eta,\gamma,\beta) \in {\cal T}\) such that
\begin{equation}\label{up_chr_case1}
\frac{1+0.5\xi}{(1-\eta-\theta_{22})(1-\beta)} \leq \frac{2(1+\xi)}{1-2\theta_2}
\end{equation}
and
\begin{equation}\label{up_chr_case1_exp}
(1-\beta)(2\eta-2\gamma)+\theta_2 \geq 1-\theta_2-\zeta.
\end{equation}
This proves the upper bound in~(\ref{chr_bd_casei}) in Theorem~\ref{chr_cor}.\\
\emph{Proof of~\((f1)\)}: We recall that~\(\theta_{22} = \frac{\theta_2}{1-\beta}\) and we have the constraint that~\(\beta > \theta_2\) and~\(\gamma > 0.\) Since
\begin{equation}\label{eta_bd}
\inf_{\beta > \theta_2,\gamma > 0} \left(\gamma + \frac{1-\theta_{22}}{2}\right) = \frac{1}{2}\left(1-\frac{\theta_2}{1-\theta_2}\right) = \frac{1-2\theta_2}{2(1-\theta_2)}
\end{equation}
we have from~(\ref{chr_condi}) that the least possible value for~\(\eta\) is~\(\frac{1-2\theta_2}{2(1-\theta_2)};\) i.e., if~\((\eta,\gamma,\beta) \in {\cal T},\) then~\(\eta \geq \frac{1-2\theta_2}{2(1-\theta_2)}\) and~\(\beta > \theta_2.\)

Fix~\(\delta > 0\) small and fix~\(\theta_2 < \beta < \theta_2 +\delta\) and~\(\frac{1-2\theta_2}{2(1-\theta_2)} < \eta < \frac{1-2\theta_2}{2(1-\theta_2)} + \delta.\) From~(\ref{eta_bd}), we have that if~\(\delta,\gamma > 0\) are sufficiently small, then~\((\eta,\gamma,\beta) \in {\cal T}.\)
Also, we have that~\((1-\eta-\theta_{22})(1-\beta) \geq f(\theta_2,\delta),\)
where~\[f(\theta_2,\delta) = \left(1-\frac{1-2\theta_2}{2(1-\theta_2)} - \delta -\frac{\theta_2}{1-\theta_2-\delta}\right)\left(1-\theta_2-\delta\right).\] Since~\(f(\theta_2,0) = \frac{1-2\theta_2}{2},\) we fix~\(\xi > 0\) and choose~\(\delta  = \delta(\xi) > 0\) smaller if necessary so that
\[(1-\theta_{22} - \eta)(1-\beta)\geq  \frac{1-2\theta_2}{2} \frac{1+0.5\xi}{1+\xi}.\] This proves~(\ref{up_chr_case1}).


For~(\ref{up_chr_case1_exp}), we proceed analogously and use the fact that~\(\eta \geq \frac{1-2\theta_2}{2(1-\theta_2)}\) to obtain
\[(1-\beta)(2\eta-2\gamma) + \theta_2 \geq (1-\theta_2-\delta)\left(\frac{1-2\theta_2}{1-\theta_2} - 2\gamma\right) + \theta_2 \geq 1-\theta_2 -\zeta,\] provided~\(\delta,\gamma > 0\) are small.~\(\qed\)







\emph{Proof of~\((ii)\)}: Here~\(r_n = p\) for some~\(p \in (0,1)\) and for all~\(n \geq 2.\) As in the proof of~\((i)\) above, we identify cliques in subsets of the random graph~\(G(n,1-r_n).\) Fix~\(\beta > 0\) to be determined later and set~\(m = n^{1-\beta},p_m = 1-p\) and apply Theorem~\ref{clq_cor}, case~\((ii)\) for the random graph~\(G(m,1-p).\)  We then have~\(\alpha_1 = \alpha_2 = 0,\) where~\(\alpha_1\) and~\(\alpha_2\) are as defined in~(\ref{alp1_def}) and~(\ref{alp2_def}), respectively. Fixing
\begin{equation}\label{chr_condii}
\frac{1+\gamma}{2} < \eta < 1,
\end{equation}
we have  from the proof of lower bound of~(\ref{clq_bd_case2}), that there is a positive integer~\(N= N(\eta,\gamma) \geq 1\) so that
\begin{eqnarray}
\mathbb{P}\left(\omega(G(m,1-p)) \leq L\right) \leq 3\exp\left(-m^{2\eta - 2\gamma}\right) \label{clq_bd_case2m}
\end{eqnarray}
for all~\(m \geq N_3,\) where
\begin{equation}\label{ldef2}
L = (1-\eta)\frac{\log{m}}{\log\left(\frac{1}{1-p}\right)} = (1-\eta)(1-\beta)\frac{\log{n}}{\log\left(\frac{1}{1-p}\right)}.
\end{equation}
The final estimate above follows using~\(m = n^{1-\beta}.\)

As in case~\((i),\) let~\(F_n\) denote the event that every set of~\(m\) vertices in the random graph~\(G(n,1-r_n)\) contains an~\(L-\)clique. Analogous to~(\ref{fn_est}), we then have
\begin{equation}\label{fn_est2}
\mathbb{P}(F_n^c) \leq {n \choose m}3\exp\left(-m^{2\eta - 2\gamma}\right) \leq n^{m} 3\exp\left(-m^{2\eta - 2\gamma}\right) = 3e^{-B},
\end{equation}
where
\begin{equation}\label{bestii}
B = m^{2\eta - 2\gamma} - m\log{n} = m^{2\eta - 2\gamma} - \frac{m}{1-\beta}\log{m}.
\end{equation}
The final estimate follows since~\(m = n^{1-\beta}.\) From the choices of~\(\eta\) and~\(\gamma\) in~(\ref{chr_condii}), we have that~\(2\eta - 2\gamma > 1\) and so
\begin{equation}\label{best2ii}
B \geq \frac{1}{2}m^{2\eta-2\gamma} = \frac{1}{2}n^{(1-\beta)(2\eta-2\gamma)}
\end{equation}
for all~\(n \geq N_2.\) Here~\(N_2 = N_2(\eta,\gamma,\beta) \geq 1\) is a constant and the final equality follows from the definition of~\(m = n^{1-\beta}.\) 

If the event~\(F_n\) occurs, then using property~\((d3),\) we have that
\begin{equation}\label{chi_uper_caseii_est}
\chi(G(n,r_n)) \leq \frac{n-m}{L} + m  \leq \frac{n}{L}+m = \frac{\log\left(\frac{1}{1-p}\right)}{(1-\eta)(1-\beta)}\frac{n}{\log{n}} + n^{1-\beta}
\end{equation}
Fixing~\(\beta > 0\) and~\(\xi > 0,\) we have that the final term in~(\ref{chi_uper_casei_est}) is at most~\[\left(\frac{1+0.5\xi}{(1-\eta)(1-\beta)}\right)\frac{n\log\left(\frac{1}{1-p}\right)}{\log{n}}\] for all~\(n  \geq N_3.\) Here~\(N_3 = N_3(\eta,\beta,\xi) \geq 1\) is a constant.
Summarizing, we have from~(\ref{bestii}),~(\ref{best2ii}) and~(\ref{fn_est2}) that
\begin{eqnarray}
&&\mathbb{P}\left(\chi(G(n,r_n)) \leq \left(\frac{1+0.5\xi}{(1-\eta)(1-\beta)}\right)\frac{n\log\left(\frac{1}{1-p}\right)}{\log{n}}\right) \nonumber\\
&&\;\;\;\;\;\geq 1 - 3\exp\left(-\frac{1}{2}n^{(1-\beta)(2\eta-2\gamma)}\right)\label{chi_est222}
\end{eqnarray}
Analogous to property~\((f1)\) above, we have the following property.\\
\((f2)\) Let \[{\cal T} = \{(\eta,\gamma,\beta) : \beta > 0 \text{ and~(\ref{chr_condii}) is satisfied}\}\] and fix~\(\zeta > 0.\)
There exists~\((\eta,\gamma,\beta) \in {\cal T}\) so that
\begin{equation}\label{up_chr_case2}
\frac{1+0.5\xi}{(1-\eta)(1-\beta)} \leq 2(1+\xi)
\end{equation}
and
\begin{equation}\label{up_chr_case2_exp}
(1-\beta)(2\eta-2\gamma) \geq 1 - \zeta.
\end{equation}
Substituting the above into~(\ref{chi_est222}), we obtain the upper bound in~(\ref{chr_bd_caseii}) in Theorem~\ref{chr_cor}.\\
\emph{Proof of~\((f2)\)}: From~(\ref{chr_condii}), we have that the minimum possible value for~\(\eta\) is~\(\frac{1}{2}.\) Choosing~\(\eta  >\frac{1}{2}\) sufficiently close to~\(\frac{1}{2}\) and~\(\beta > 0\) small, both~(\ref{up_chr_case2}) and~(\ref{up_chr_case2_exp}) are satisfied.\(\qed\)



\emph{Proof of~\((iii)\)}: Here~\(r_n = 1-\frac{1}{n^{\theta_1}}\) for some~\(\theta_1 > 0.\) As before, we identify cliques in subsets of the random graph~\(G(n,1-r_n).\) Fix~\(\beta > 0\) to be determined later and set~\(m = \beta n\) and apply Theorem~\ref{clq_cor}, case~\((iii)\) for the random graph~\(G(m,p_m),\) where
\begin{equation}\label{pm_def3}
p_m = 1-\frac{\beta^{\theta_1}}{m^{\theta_1}} = 1-r_n.
\end{equation}
We then have~\(\alpha_1 = \theta_1\) and~\(\alpha_2 = 0,\) where~\(\alpha_1\) and~\(\alpha_2\) are as defined in~(\ref{alp1_def}) and~(\ref{alp2_def}), respectively. Let~\(\eta,\gamma >0\) be such that
\begin{equation}\label{chr_condiii}
\frac{1+\theta_{1}}{2}  + \gamma < \eta <1.
\end{equation}

From the proof of lower bound of~(\ref{clq_bd_case1}), we have that
\begin{eqnarray}
\mathbb{P}\left(\omega(G(m,p_m)) \leq L\right) \leq 3\exp\left(-m^{2\eta - 2\gamma-\theta_{1}}\right) \label{clq_bd_case3m}
\end{eqnarray}
for all~\(m\) large, where
\begin{equation}\label{ldef3}
L = \frac{2(1-\eta)}{\theta_1}.
\end{equation}

As in cases~\((i)-(ii),\) let~\(F_n\) denote the event that every set of~\(m\) vertices in the random graph~\(G(n,1-r_n)\) contains an~\(L-\)clique. Using~\({n \choose k} \leq \left(\frac{ne}{k}\right)^{k}\) for integers~\(1 \leq k \leq n,\) we have
\begin{equation}\label{fn_est3}
\mathbb{P}(F_n^c) \leq {n \choose m}3\exp\left(-m^{2\eta - 2\gamma-\theta_{1}}\right) \leq \left(\frac{ne}{m}\right)^{m} 3\exp\left(-m^{2\eta - 2\gamma-\theta_{1}}\right) = 3e^{-B},
\end{equation}
where
\begin{equation}\label{bestiii}
B = m^{2\eta - 2\gamma-\theta_{1}} - m\log\left(\frac{e}{\beta}\right).
\end{equation}
The final estimate in~(\ref{fn_est3}) follows since~\(m = \beta n.\) From the choices of~\(\eta\) and~\(\gamma\) in~(\ref{chr_condiii}), we have that~\(2\eta - 2\gamma-\theta_{1} > 1\) and so
\begin{equation}\label{best3iii}
B \geq \frac{1}{2}m^{2\eta-2\gamma-\theta_{1}} = \frac{1}{2}\left(\beta n\right)^{2\eta-2\gamma-\theta_1}
\end{equation}
for all~\(n \geq N_2.\) Here~\(N_2 = N_2(\eta,\gamma,\beta) \geq 1\) is a constant and the final equality follows from the definition of~\(m = \beta n.\)

If the event~\(F_n\) occurs, then using property~\((d3),\) we have that
\begin{equation}\label{chi_uper_caseiii_est}
\chi(G(n,r_n)) \leq \frac{n}{L} + 2m  \leq \frac{\theta_1}{2(1-\eta)}n + \beta n + 1.
\end{equation}
Summarizing, we have from~(\ref{best3iii}) and~(\ref{fn_est3}) that
\begin{equation}\label{chi_est333}
\mathbb{P}\left(\chi(G(n,r_n)) \leq \frac{\theta_1}{2(1-\eta)}n + \beta n+1\right) \geq 1-\exp\left(-\frac{1}{2}\left(\beta n\right)^{2\eta-2\gamma-\theta_1}\right)
\end{equation}
for all~\(n\) large.
We have the following property.\\
\((f3)\) Fix~\(\xi,\zeta > 0\) and let \[{\cal T} = \{(\eta,\gamma,\beta) : \text{The condition~(\ref{chr_condiii}) is satisfied}\}.\]
There exists~\((\eta,\gamma,\beta) \in {\cal T}\) such that
\begin{equation}\label{up_chr_case3}
\frac{\theta_1}{2(1-\eta)} +\beta \leq \frac{\theta_1}{1-\theta_1} (1+\xi)
\end{equation}
and
\begin{equation}\label{up_chr_case3_exp}
\frac{1}{2}\left(\beta n\right)^{2\eta-2\gamma-\theta_1} \geq n^{1-\theta_1-\zeta}
\end{equation}
for all~\(n\) large. Using the above estimates in~(\ref{chi_uper_caseiii_est}), we obtain the upper bound in~(\ref{chr_bd_caseiii}) in Theorem~\ref{chr_cor}.\\
\emph{Proof of~\((f3)\)}: From~(\ref{chr_condiii}), we have that the least possible value for~\(\eta\) is~\(\frac{1+\theta_1}{2}.\) Fix~\(\zeta,\xi > 0.\)  Choosing~\(\gamma,\beta > 0\) sufficiently small and~\(\eta > \frac{1+\theta_1}{2}\) sufficiently close of~\(\frac{1+\theta_1}{2},\) we get that
\[\frac{\theta_1}{2(1-\eta)} +\beta  \leq \frac{\theta_1}{1-\theta_1} (1+\xi)\] and~\(2\eta - 2\gamma - \theta_1 >1-\zeta.\) This proves~(\ref{up_chr_case3}) and~(\ref{up_chr_case3_exp}).~\(\qed\)

\subsection*{Acknowledgement}
I thank Professors Rahul Roy and Federico Camia for crucial comments and for my fellowships.

\bibliographystyle{plain}

\end{document}